\documentclass[11pt,a4paper]{article}

\usepackage[
  a4paper,
  left=2cm,
  right=2cm,
  top=2cm,
  bottom=2cm,
  footskip=1cm
]{geometry}

\usepackage[T1]{fontenc}
\usepackage{lmodern}
\usepackage{microtype}
\usepackage{amsmath,amssymb,amsthm}
\usepackage{mathtools}
\usepackage{aliascnt}

\usepackage{graphicx}
\usepackage{float}

\usepackage{xcolor}
\definecolor{interactionblue}{HTML}{216A9D}

\usepackage{cite}
\usepackage[shortlabels]{enumitem}

\usepackage{booktabs,array,tabularx,longtable}
\newcolumntype{Y}{>{\raggedright\arraybackslash}X}

\usepackage{caption}
\usepackage{tikz}
\usetikzlibrary{positioning,calc,arrows.meta}

\usepackage[hidelinks]{hyperref}
\usepackage[nameinlink,capitalise,noabbrev]{cleveref}

\theoremstyle{plain}

\newtheorem{theorem}{Theorem}[section]

\newaliascnt{lemma}{theorem}
\newtheorem{lemma}[lemma]{Lemma}
\aliascntresetthe{lemma}

\newaliascnt{proposition}{theorem}

\aliascntresetthe{proposition}

\newaliascnt{corollary}{theorem}
\newtheorem{corollary}[corollary]{Corollary}
\aliascntresetthe{corollary}

\newaliascnt{conjecture}{theorem}
\newtheorem{conjecture}[conjecture]{Conjecture}
\aliascntresetthe{conjecture}

\newaliascnt{claim}{theorem}
\newtheorem{claim}[claim]{Claim}
\aliascntresetthe{claim}

\newaliascnt{problem}{theorem}

\aliascntresetthe{problem}

\newaliascnt{question}{theorem}

\aliascntresetthe{question}

\theoremstyle{definition}

\newaliascnt{definition}{theorem}

\aliascntresetthe{definition}

\newaliascnt{example}{theorem}

\aliascntresetthe{example}

\newaliascnt{configuration}{theorem}

\aliascntresetthe{configuration}

\theoremstyle{remark}

\newaliascnt{remark}{theorem}
\newtheorem{remark}[remark]{Remark}
\aliascntresetthe{remark}

\crefname{theorem}{Theorem}{Theorems}
\Crefname{theorem}{Theorem}{Theorems}

\crefname{lemma}{Lemma}{Lemmas}
\Crefname{lemma}{Lemma}{Lemmas}

\crefname{proposition}{Proposition}{Propositions}
\Crefname{proposition}{Proposition}{Propositions}

\crefname{corollary}{Corollary}{Corollaries}
\Crefname{corollary}{Corollary}{Corollaries}

\crefname{conjecture}{Conjecture}{Conjectures}
\Crefname{conjecture}{Conjecture}{Conjectures}

\crefname{claim}{Claim}{Claims}
\Crefname{claim}{Claim}{Claims}

\crefname{problem}{Problem}{Problems}
\Crefname{problem}{Problem}{Problems}

\crefname{question}{Question}{Questions}
\Crefname{question}{Question}{Questions}

\crefname{definition}{Definition}{Definitions}
\Crefname{definition}{Definition}{Definitions}

\crefname{example}{Example}{Examples}
\Crefname{example}{Example}{Examples}

\crefname{configuration}{Configuration}{Configurations}
\Crefname{configuration}{Configuration}{Configurations}

\crefname{remark}{Remark}{Remarks}
\Crefname{remark}{Remark}{Remarks}

\crefname{figure}{Fig.}{Figs.}
\Crefname{figure}{Fig.}{Figs.}

\newcommand{\Nstr}{N^{\mathrm s}}
\newcommand{\Azero}{A_0}
\newcommand{\colset}{[28]}

\DeclareMathOperator{\dist}{dist}
\newcommand{\rhoR}{\rho_R}

\date{}

\begin{document}

\baselineskip=17pt
\parindent=0.5cm

\title{\Large {\bf Strong chromatic index of graphs with Ore-degree at most nine}}
\bigskip

\author{%
Xiaoxue Hu\textsuperscript{1}, Jiangxu Kong\textsuperscript{2}\thanks{Corresponding author.
E-mail: \texttt{kjx@hznu.edu.cn}.} and Yu Xue\textsuperscript{1}
\\[0.5em]
\small\textsuperscript{1}School of Science,
Zhejiang University of Science and Technology,
Hangzhou 310023, China\\
\small\textsuperscript{2}School of Mathematics, Hangzhou Normal University, Hangzhou 311121, China\\
}

\date{}
\maketitle

%%%%%%%%%%%%%%%%%%%%%%%%%%%%%%%%%%%%%%%%%%%%%%%%%%%%%%%%%%%%%%%%%%%%%%%
%
\begin{abstract}

\baselineskip=14 pt
The strong chromatic index $\chi'_s(G)$ of a graph $G$ is the
minimum number of colors in an edge-coloring in which distinct edges
at distance at most two in the line graph receive different colors.
Let $\theta(G)=\max_{xy\in E(G)}(d_G(x)+d_G(y))$ denote the
maximum edge weight, also called the Ore-degree, of $G$.
We prove that every finite loopless multigraph $G$ with $\theta(G)\le9$
satisfies $\chi'_s(G)\le28$.
For simple graphs, this improves the direct greedy bound from $32$ to $28$.
In particular, it verifies the Erd\H{o}s--Ne\v{s}et\v{r}il conjecture
for simple graphs of maximum degree five whose vertices of degree five
form an independent set, using one fewer color than the conjectured
bound of $29$.
When $\theta(G)=9$, our bound is four colors above the bound of $24$
predicted by the Ore-degree conjecture of Chen, Huang, Yu, and Zhou.

\par\medskip
\noindent\textbf{Keywords:} strong edge-coloring; strong chromatic index; maximum edge weight; Ore-degree; multigraphs.

\par\noindent\textbf{2020 Mathematics Subject Classification:} 05C15, 05C35.

\end{abstract}
%
%%%%%%%%%%%%%%%%%%%%%%%%%%%%%%%%%%%%%%%%%%%%%%%%%%%%%%%%%%%%%%%%%%%%%%%

\baselineskip=17pt

\section{Introduction}\label{sec:introduction}

All graphs are finite, undirected, and loopless; parallel edges
are allowed unless stated otherwise.
For a vertex $v$ of a graph $G$, let $d_G(v)$ denote its degree, with parallel edges counted with multiplicity, and let $\Delta(G)$ denote the maximum degree of $G$.

The \emph{line graph} $L(G)$ of a graph $G$ is the simple graph with vertex set $E(G)$ in which two distinct vertices
are adjacent if and only if the corresponding edges of $G$ share an endpoint. For a positive integer $k$, a \emph{strong $k$-edge-coloring} of a graph $G$ is a map $\varphi:E(G)\to\{1,\ldots,k\}$ such that $\varphi(e)\ne\varphi(f)$ whenever distinct edges $e$ and $f$ have distance at most two in the line graph $L(G)$. Equivalently, $\varphi$ is a proper vertex-coloring of $L(G)^2$. Here $L(G)^2$ is the square of $L(G)$, obtained by joining distinct vertices whose distance in $L(G)$ is at most two. When $G$ is simple, this is also equivalent to requiring each color class to be an induced matching. The \emph{strong chromatic index} of $G$,
denoted by $\chi'_s(G)$, is the minimum number of colors needed for a strong edge-coloring of $G$. Fouquet and Jolivet~\cite{FJ} introduced the notion of strong edge-coloring.

For a graph of maximum degree $\Delta$, a greedy coloring gives $\chi'_s(G)\le2\Delta^2-2\Delta+1$. Erd\H{o}s and Ne\v{s}et\v{r}il~\cite{EN} proposed the following conjecture.

\begin{conjecture}
\label{conj:erdos-nesetril}
Let $G$ be a finite simple graph with maximum degree $\Delta$. Then
\[
  \chi'_s(G)\le
  \begin{cases}
    \dfrac54\Delta^2,
      & \Delta\text{ is even},\\[2mm]
    \dfrac14(5\Delta^2-2\Delta+1),
      & \Delta\text{ is odd}.
  \end{cases}
\]
\end{conjecture}

\Cref{conj:erdos-nesetril} holds for $\Delta\le3$. The case $\Delta=3$ was proved independently by Andersen~\cite{Andersen} and by Hor\'{a}k, He, and Trotter~\cite{Horak}. Huang, Santana, and Yu~\cite{Huang} proved $\chi'_s(G)\le21$ for simple graphs with maximum degree at most $4$, and Zang~\cite{Zang} proved $\chi'_s(G)\le37$ for simple graphs with maximum degree at most $5$. For sufficiently large $\Delta$, Molloy and Reed~\cite{MolloyReed} proved $\chi'_s(G)\le1.998\Delta^2$. This coefficient was subsequently reduced to $1.93$ by Bruhn and Joos~\cite{BruhnJoos}, to $1.835$ by Bonamy, Perrett, and Postle~\cite{Bonamy}, and to $1.772$ by Hurley, de Joannis de Verclos, and Kang~\cite{HurleyEtAl}. In a 2026 preprint, Davey, Hurley, de Joannis de Verclos, Kang, and Volec~\cite{DaveyEtAl2026} proved $\chi'_s(G)\le1.73\Delta^2$ when $\Delta$ is sufficiently large.

We study strong edge-coloring in terms of the degrees of adjacent vertices. The \emph{weight} of an edge $xy$ is $d_G(x)+d_G(y)$, and the \emph{maximum edge weight}, also called the \emph{maximum Ore-degree}, is $\theta(G)=\max_{xy\in E(G)}\bigl(d_G(x)+d_G(y)\bigr)$. For an edgeless graph, set $\theta(G)=\chi'_s(G)=0$. If $G$ is simple and $E(G)\ne\varnothing$, then the vertex of $L(G)$ corresponding to $xy$ has degree $d_G(x)+d_G(y)-2$, and hence $\theta(G)=\Delta(L(G))+2$. The parameter $\theta(G)$ records the largest degree sum across an edge. Although $\theta(G)\le2\Delta(G)$, the inequality can be strict when vertices of large degree are nonadjacent. Bounds in terms of $\theta(G)$ can therefore provide more information than bounds based only on $\Delta(G)$.

Chen, Huang, Yu, and Zhou~\cite{Chen2020} proposed the following Ore-degree analogue of \cref{conj:erdos-nesetril}.

\begin{conjecture}
\label{conj:ore-degree}
Let $G$ be a finite simple graph with maximum edge weight $\theta=\theta(G)\ge5$, and set $q=\lceil\theta/4\rceil$. Then $\chi'_s(G)\le B(\theta)$, where
\[
  B(\theta)=
  \begin{cases}
    5q^2-8q+3, & \theta\equiv1\pmod4,\\
    5q^2-6q+2, & \theta\equiv2\pmod4,\\
    5q^2-4q+1, & \theta\equiv3\pmod4,\\
    5q^2,      & \theta\equiv0\pmod4.
  \end{cases}
\]
\end{conjecture}

Wu and Lin~\cite{WuLin} proved that every simple graph $G$ with $\theta(G)\le5$ satisfies $\chi'_s(G)\le7$. More precisely, the bound is $6$ if no component of $G$ is isomorphic to the graph obtained from a $5$-cycle by adding a vertex adjacent to two nonadjacent vertices of the cycle; the exceptional graph has strong chromatic index $7$. Nakprasit and Nakprasit~\cite{Nakprasit} proved that every loopless multigraph $G$ with $\theta(G)\le6$ satisfies $\chi'_s(G)\le10$, with the stronger bound $9$ when $G$ is bipartite. Chen, Huang, Yu, and Zhou~\cite{Chen2020} also proved the bound $10$ for simple graphs with maximum edge weight at most $6$, and they proved $\chi'_s(G)\le15$ when $\theta(G)\le7$, compared with the conjectured bound $B(7)=13$. For simple graphs with $\theta(G)\le8$, Chen, Chen, Zhao, and Zhou~\cite{ChenEight} proved $\chi'_s(G)\le21$, one more than the conjectured bound $B(8)=20$.

We prove the following bound for graphs with maximum edge weight at most nine, allowing parallel edges.

\begin{theorem}\label{thm:main}
Every finite loopless multigraph $G$ with $\theta(G)\le9$ satisfies $\chi'_s(G)\le28$.
\end{theorem}

For simple graphs with $\theta(G)=9$, \cref{conj:ore-degree} predicts $\chi'_s(G)\le B(9)=24$. A direct greedy argument gives the bound $32$: by \cref{lem:strong-neighborhood}, every edge of a simple graph with $\theta(G)\le9$ has at most $31$ other edges at distance at most two in its line graph. \Cref{thm:main} improves this upper bound to $28$, thereby halving the gap between the direct greedy bound and the conjectured bound. Moreover, if $\Delta(G)\le5$ and no two vertices of degree $5$ are adjacent, then every edge of $G$ has weight at most nine. Thus, \cref{thm:main} yields the following corollary.

\begin{corollary}\label{cor:independent-five-vertices}
Let $G$ be a finite loopless multigraph with $\Delta(G)\le5$. If no two vertices of degree $5$ are adjacent, then $\chi'_s(G)\le28$.
\end{corollary}

For simple graphs of maximum degree five whose vertices of degree five form an independent set, this gives one fewer color than the bound $29$ in \cref{conj:erdos-nesetril}.

We conclude this section with an outline of the proof.
We choose a counterexample $G$ minimizing $|V(G)|+|E(G)|$.
Deletion and extension arguments show that $G$ is simple, every
vertex has degree $4$ or $5$, $G$ has no edge-cut of size at most
three, and $G$ has girth at least six.
The bound of Huang, Santana, and Yu~\cite{Huang}
then ensures that $G$ has a vertex of degree five. This is the
only external bound on the strong chromatic index used in our proof.
After precoloring some edges near a fixed $5$-vertex, we choose
a longest sequence of distinct edges
satisfying specified conditions.
This sequence determines a partition of the vertex set into
sets $L,M,R$ with no edge joining $L$ to $R$.
We then discard the initial precoloring and,
by minimality, independently color the induced subgraphs
$G[L]$ and $G[R]$ with $28$ colors.
Structural lemmas allow us to choose a
vertex of degree four in $M$ whose incident
edges will be colored last.
We permute the colors in $G[L]$ and $G[R]$ and extend the resulting
partial coloring,
using uncolored conflicting edges and repetitions
among the colors on colored conflicting edges to bound the number
of forbidden colors at each step.
The proofs excluding short cycles are given in
\cref{app:short-cycles}.

\section{Notation and preliminary lemmas}\label{sec:prelim}

For a vertex $v\in V(G)$, let $N_G(v)$ denote its set of neighbors.
For $V(G)\ne\varnothing$, let
$\delta(G)=\min_{v\in V(G)}d_G(v)$ denote the minimum degree of $G$.
A vertex of degree $k$ is called a $k$-vertex.
An $(a,b)$-edge has endpoints of degrees $a$ and $b$.

For $X\subseteq V(G)$, let $G[X]$ denote the
subgraph of $G$ induced by $X$.
For disjoint vertex sets $X,Y\subseteq V(G)$,
let $E(X,Y)$ denote the set of edges of $G$
with one endpoint in each set.
For any graph $Q$ and $v\in V(Q)$, let $E_Q(v)$ denote the set
of edges of $Q$ incident with $v$.
For a positive integer $m$, write $[m]=\{1,\ldots,m\}$.

For edges $e,f\in E(G)$, their \emph{edge-distance} is
$\dist_G(e,f)=\dist_{L(G)}(e,f)$, the distance
between the corresponding vertices of the line graph.
We set this distance to $\infty$ when $e$ and $f$ lie in
different components of $G$.
The \emph{strong neighborhood} of an edge $e$ is
$\Nstr_G(e)=\{f\in E(G)\setminus\{e\}:\dist_G(e,f)\le2\}$.
Edges $e$ and $f$ \emph{conflict} if $f\in\Nstr_G(e)$;
equivalently, we say that $e$ \emph{sees} $f$.

\begin{lemma}\label{lem:strong-neighborhood}
Let $xy$ be an edge of a simple graph $G$
with $\theta(G)\le9$. Then
\begin{equation}\label{eq:local-bound}
   |\Nstr_G(xy)|
   \le 
   (d_G(x)-1)(9-d_G(x))
   +(d_G(y)-1)(9-d_G(y))
   \le31.
\end{equation}
If $d_G(x)=d_G(y)=4$, then
$|\Nstr_G(xy)|\le30$.
\end{lemma}

\begin{proof}
There are $d_G(x)+d_G(y)-2$ edges other than
$xy$ incident with $x$ or $y$.
For each $x'\in N_G(x)\setminus\{y\}$, the inequality
$d_G(x')\le9-d_G(x)$ implies that at most
$d_G(x')-1\le8-d_G(x)$ additional edges
are reached through $xx'$.
The $x$-side therefore contributes at most
$(d_G(x)-1)(8-d_G(x))$
edges at distance two.
The same argument gives at most
$(d_G(y)-1)(8-d_G(y))$ such edges
on the $y$-side.
Any overlap only decreases the number of distinct edges
in $\Nstr_G(xy)$, which proves
the first inequality in \eqref{eq:local-bound}.
Moreover,
\(
{
 (d_G(x)-1)(9-d_G(x))+(d_G(y)-1)(9-d_G(y))=32-(d_G(x)-5)^2-(d_G(y)-5)^2.
}
\)
Since $d_G(x)$ and $d_G(y)$ are positive
integers with $d_G(x)+d_G(y)\le9$, we have
$(d_G(x)-5)^2+(d_G(y)-5)^2\ge1$.
Thus the expression is at most $31$, with equality if and only if
$(d_G(x),d_G(y))\in\{(4,5),(5,4)\}$.
Substituting $d_G(x)=d_G(y)=4$
into \eqref{eq:local-bound} gives
$|\Nstr_G(xy)|\le15+15=30$.
\end{proof}

{
A partial strong edge-coloring of $G$ assigns colors to a subset
of $E(G)$ so that any two colored edges that conflict in $G$
receive different colors.
All conflicts are measured in $G$.
For $U\subseteq V(G)$, every strong edge-coloring of $G[U]$
is a partial strong edge-coloring of $G$: any edge of $G$
joining endpoints of two edges of $G[U]$ also belongs to $G[U]$.
This need not hold for a subgraph obtained by deleting edges.
An extension preserves every color already assigned.
}

Let $\varphi$ be a partial strong edge-coloring with color set
$\colset=\{1,\ldots,28\}$.
For an uncolored edge $e$, define
$A_\varphi(e)=\colset\setminus
\{\varphi(f):f\in\Nstr_G(e)\text{ and }f\text{ is colored}\}$.
These are the available colors for $e$.
When the partial coloring is fixed, write $A(e)$ for $A_\varphi(e)$.

\begin{lemma}\label{lem:hall}
Let $X$ be a set of uncolored edges.
If $\left|\bigcup_{e\in Y}A(e)\right|\ge |Y|$ for every
$Y\subseteq X$, then the edges in $X$ can be assigned pairwise
distinct available colors, extending the current
partial strong edge-coloring.
In particular, this holds if the edges of $X$ can be indexed
as $e_1,\ldots,e_k$ so that $|A(e_i)|\ge i$ for $1\le i\le k$.
\end{lemma}

\begin{proof}
The first assertion follows from Hall's theorem~\cite{Hall}.
For the second, let $Y\subseteq X$ be nonempty and put
$j=\max\{i:e_i\in Y\}$.
Then $|Y|\le j$, and hence
$\left|\bigcup_{e\in Y}A(e)\right|
 \ge |A(e_j)|\ge j\ge |Y|$.
Thus Hall's condition holds.
\end{proof}

The following grouped form permits pairwise nonconflicting edges
to receive a common color.

\begin{lemma}\label{lem:grouped-hall}
Let $X_1,\ldots,X_m$ be pairwise disjoint
nonempty sets of uncolored edges, and suppose
that the edges in each $X_i$ are pairwise nonconflicting.
Put $B_i=\bigcap_{e\in X_i}A(e)$.
If $\left|\bigcup_{i\in I}B_i\right|\ge |I|$ for every
$I\subseteq[m]$, then there are distinct colors $c_i\in B_i$
such that every edge in $X_i$ can be colored with $c_i$, extending the current partial strong edge-coloring.
\end{lemma}

\begin{proof}
Apply Hall's theorem to the family $B_1,\ldots,B_m$.
Edges within a set $X_i$ do not conflict and hence may receive
their common representative $c_i$; distinct sets receive distinct colors.
Since $c_i\in B_i$, this color is available
for every edge in $X_i$.
\end{proof}

In the next three lemmas, all counts refer
to the current partial strong $28$-edge-coloring of $G$
and are evaluated immediately before the edge under consideration
is colored.

Suppose that, under a partial
strong $28$-edge-coloring $\varphi$,
an uncolored edge $e$ conflicts with $m$ colored edges
using $q$ distinct colors.
Put $r_\varphi(e)=m-q$.
Thus $r_\varphi(e)$ is the difference between the number
of colored edges in $\Nstr_G(e)$ and the number of colors they use.

\begin{lemma}\label{lem:certificate}
Let $e$ be uncolored in a partial strong $28$-edge-coloring
$\varphi$ of $G$, and suppose that
$|\Nstr_G(e)|\le B(e)$.
When $e$ is to be colored, let $t(e)$ be a lower bound
on the number of uncolored edges in $\Nstr_G(e)$,
and let $r(e)$ be a lower bound on $r_\varphi(e)$.
If $B(e)-t(e)-r(e)\le27$, then $A(e)\ne\varnothing$.
\end{lemma}

\begin{proof}
At most $B(e)-t(e)$ edges in $\Nstr_G(e)$ are colored.
These edges use at most
$B(e)-t(e)-r(e)$ distinct colors. Hence
$|A(e)|\ge28-B(e)+t(e)+r(e)\ge1.$
\end{proof}

\begin{lemma}\label{lem:four-savings}
Let $e$ be uncolored in a partial strong
$28$-edge-coloring $\varphi$ of $G$, and let $t(e),r(e)$
be lower bounds as in \cref{lem:certificate}.
If $|\Nstr_G(e)|\le31$ and $t(e)+r(e)\ge4$,
then $A(e)\ne\varnothing$.
\end{lemma}

\begin{proof}
Apply \cref{lem:certificate} with $B(e)=31$.
\end{proof}

\begin{lemma}\label{lem:three-savings}
Let $G$ be a simple graph with $\theta(G)\le9$,
and let $e=xy$ be uncolored in a partial strong
$28$-edge-coloring $\varphi$ of $G$.
Suppose that $d_G(x)=d_G(y)=4$, and let $t(e),r(e)$
be lower bounds as in \cref{lem:certificate}.
If $t(e)+r(e)\ge3$, then $A(e)\ne\varnothing$.
\end{lemma}

\begin{proof}
By \cref{lem:strong-neighborhood}, $|\Nstr_G(e)|\le30$.
Apply \cref{lem:certificate} with $B(e)=30$.
\end{proof}

Let $\varphi:E(G)\to\colset$ be an edge-coloring
of $G$, and let $\pi:\colset\to\colset$ be a permutation
of the color set.
We denote by $\pi\circ\varphi$ the edge-coloring of $G$ defined by
$(\pi\circ\varphi)(e)=\pi(\varphi(e))$ for every $e\in E(G)$.
Thus, every edge of color $c$ under $\varphi$ receives
color $\pi(c)$ under $\pi\circ\varphi$.
The same definition applies to a partial coloring by restricting
the domain to its colored edges.

\begin{lemma}\label{lem:palette-permutation}
If $\varphi$ is a strong $28$-edge-coloring
of $G$ and $\pi$ is a permutation of $\colset$,
then $\pi\circ\varphi$ is also a strong $28$-edge-coloring
of $G$.
\end{lemma}

\begin{proof}
Since $\pi$ is a bijection, two edges have the same color under
$\pi\circ\varphi$ if and only if they have the same color
under $\varphi$.
Thus $\pi\circ\varphi$ is a strong edge-coloring.
\end{proof}

We use the following convention when constructing auxiliary graphs.

\begin{remark}\label{rem:auxiliary}
{
In the proof of \cref{thm:main}, we choose a counterexample $G$
minimizing $|V(G)|+|E(G)|$.
Whenever we apply minimality to an auxiliary graph $Q$,
we check that $Q$ is loopless, $\theta(Q)\le9$, and
$|V(Q)|+|E(Q)|<|V(G)|+|E(G)|$.

When transferring a coloring of $Q$ back to $G$, we retain colors
only on edges belonging to both graphs.
We record any colors of added edges that will be used later
and discard the added edges themselves.
Before extending the retained coloring, we verify that it is
a partial strong edge-coloring of $G$.

The color of an added edge may be assigned to an uncolored edge
$e$ of $G$ only after checking that it differs from the colors
on all edges in $\Nstr_G(e)$ that are colored at the time
of the assignment.
If several edges are colored simultaneously, we also check
that any two of them that conflict receive different colors.
}
\end{remark}

\section{Structure of a minimum counterexample}\label{sec:minimal}

Assume that \cref{thm:main} is false, and choose a counterexample
$G$ minimizing $|V(G)|+|E(G)|$.
By minimality, every loopless multigraph $Q$ with $\theta(Q)\le9$
and $|V(Q)|+|E(Q)|<|V(G)|+|E(G)|$
has a strong $28$-edge-coloring.
The graph $G$ is connected and has no isolated vertices, since
distinct components can be colored independently with the same
set of colors.
Since $\theta(G)\le9$, every vertex of $G$ has degree at most $8$.

\begin{claim}\label{cl:simple}
The graph $G$ is simple.
\end{claim}

\begin{proof}
Suppose that $e$ and $f$ are parallel edges
with endpoints $x$ and $y$.
By minimality, $G-e$ has a strong $28$-edge-coloring.
Since $f$ remains, deleting $e$ does not remove any conflict
between retained edges, so this is a partial strong
edge-coloring of $G$.

There are at most $d_G(x)+d_G(y)-3$ edges
other than $e$ incident with $x$ or $y$.
Using $\theta(G)\le9$ to count the additional edges reached
through their other neighbors gives
\[
{
\begin{aligned}
|\Nstr_G(e)|
&\le d_G(x)+d_G(y)-3+(d_G(x)-2)(8-d_G(x))+(d_G(y)-2)(8-d_G(y))\\
&=1+(d_G(x)-2)(9-d_G(x))
    +(d_G(y)-2)(9-d_G(y))\\
&\le25,
\end{aligned}
}
\]
since each product in the preceding expression
is at most $12$.
Thus a color is available for $e$, a contradiction.
\end{proof}

\begin{claim}\label{cl:degrees}
Every vertex of $G$ has degree four or five.
\end{claim}

\begin{proof}
Let $u$ be a vertex of degree $k$, write
$N_G(u)=\{u_1,\ldots,u_k\}$, and put $d_i=d_G(u_i)$.
By minimality, $G-u$ has a strong
$28$-edge-coloring, which is a partial strong edge-coloring
of $G$ since $G-u$ is induced.
For each $i$, the edge $uu_i$ sees at most
\begin{equation}\label{eq:star-bound}
  (d_i-1)(9-d_i)+\sum_{j\ne i}(d_j-1)
\end{equation}
colored edges.
The first term bounds the colored conflicts through $u_i$;
the sum bounds the colored edges incident with the other
neighbors of $u$.

Suppose that $k\le3$, and number the neighbors
so that $d_k=\max_{1\le i\le k}d_i$.
Since $d_i\le9-k$ for every $i$, the quantity in
\eqref{eq:star-bound} is at most
\(
  16+(k-1)(8-k)\le26.
\)
Thus $|A(uu_i)|\ge2$ for every $i$.
For $uu_k$, the same quantity is at most
\(
  (d_k-1)(11-d_k)
  =25-(d_k-6)^2
  \le25,
\)
because $k\le3$ and $d_j\le d_k$ for every $j$.
Hence $|A(uu_k)|\ge3$.
By \cref{lem:hall}, the star edges can be assigned pairwise
distinct available colors, extending the coloring to $G$,
a contradiction.

Therefore $\delta(G)\ge4$.
Since $\theta(G)\le9$, every vertex has degree at most $5$,
which proves the claim.
\end{proof}

\begin{claim}\label{cl:cut}
The graph $G$ has no edge-cut of size at most three.
\end{claim}

\begin{proof}

Suppose that $S=\{e_1,\ldots,e_t\}$ is an inclusion-minimal
edge-cut with $1\le t\le3$, and let $G_1,G_2$ be the two
components of $G-S$.

For $j\in\{1,2\}$, obtain $Q_j$ from $G$ by contracting
$G_{3-j}$ to a vertex $z_j$, deleting loops and retaining
parallel edges.
Let $e_i^{(j)}$ denote the edge corresponding to $e_i$ in $Q_j$.
Every vertex of $G_j$ has the same degree in $Q_j$ as in $G$,
while $d_{Q_j}(z_j)=t\le3$.
Thus $Q_j$ is loopless and $\theta(Q_j)\le9$.
Since $\delta(G)\ge4>t$, both components have at least two
vertices, so
$|V(Q_j)|+|E(Q_j)|<|V(G)|+|E(G)|$.
By minimality, $Q_j$ has a strong $28$-edge-coloring $\varphi_j$.

By permuting the colors in $Q_2$, we may assume that
$\varphi_1(e_i^{(1)})=\varphi_2(e_i^{(2)})=c_i$ $(i\in[t])$,
where $c_1,\ldots,c_t$ are distinct.
For $j\in\{1,2\}$, let $\mathcal{C}_j$ be the set of colors
used by $\varphi_j$ on edges of $G_j$ incident with
endpoints of edges in $S$.
Every such edge in $G_j$ conflicts with every $e_i^{(j)}$,
so $\{c_1,\ldots,c_t\}\cap
(\mathcal{C}_1\cup\mathcal{C}_2)=\varnothing$.
Since each endpoint of an edge in $S$ is incident with
at most four edges in its component,
$|\mathcal{C}_j|\le4t\le12$.
As $28-t-|\mathcal{C}_1|\ge13>|\mathcal{C}_2|$,
a further permutation of the colors in $Q_2$, fixing each
$c_i$, makes the two boundary color sets disjoint.

Retain the resulting colorings on $G_1$ and $G_2$, and
assign the distinct colors $c_i$ to the cut edges $e_i$.
Conflicts within either component are respected because
each $G_j$ is induced in $G$.
Every edge of $G_j$ conflicting with $e_i$ in $G$ also
conflicts with $e_i^{(j)}$ in $Q_j$, and hence avoids $c_i$.
Finally, conflicting edges in different components are
incident with opposite ends of a cut edge, so their colors
differ by the choice of the second permutation.
We obtain a strong $28$-edge-coloring of $G$, a contradiction.
\end{proof}

\begin{claim}\label{cl:girth}
The graph $G$ has girth at least six.
\end{claim}

\begin{proof}
By
\cref{lem:no-triangle-app,lem:no-c4-app,lem:no-c5-app},
$G$ contains no triangle, $4$-cycle, or $5$-cycle.
Since $G$ is simple, its girth is at least six.
\end{proof}

\begin{figure}[t]
\centering
\begin{tikzpicture}[
  x=1cm,
  y=1cm,
  vtx/.style={
    circle,
    fill=black,
    inner sep=1.5pt
  },
  lab/.style={
    font=\small,
    inner sep=0pt
  },
  edge/.style={
    line width=.8pt,
    line cap=round
  },
  optional/.style={
    edge,
    dashed
  }
]

% Root
\node[vtx] (u) at (0,4.5) {};
\node[lab,anchor=south]
  at ($(u)+(0,.16)$) {$u$};

% Five branches, arranged symmetrically
\foreach \name/\X in {
  v/-5.6,
  w/-2.8,
  x/0,
  y/2.8,
  z/5.6
}{
  \node[vtx] (\name) at (\X,3) {};
  \draw[edge] (u)--(\name);

  % First-level labels
  \ifdim\X pt<0pt
    \node[lab,anchor=east]
      at ($(\name)+(-.16,0)$) {$\name$};
  \else
    \node[lab,anchor=west]
      at ($(\name)+(.16,0)$) {$\name$};
  \fi

  % Three second-level vertices in each branch
  \foreach \i/\dx in {
    1/-.85,
    2/0,
    3/.85
  }{
    \node[vtx] (\name\i)
      at ($(\name)+(\dx,-1.45)$) {};

    \draw[edge] (\name)--(\name\i);

    % Put labels on the outward side
    \pgfmathsetmacro{\childX}{\X+\dx}

    \ifdim\childX pt<0pt
      \def\dir{-1}
      \node[lab,anchor=east]
        at ($(\name\i)+(-.16,0)$)
        {$\name_{\i}$};
    \else
      \def\dir{1}
      \node[lab,anchor=west]
        at ($(\name\i)+(.16,0)$)
        {$\name_{\i}$};
    \fi

    % Three solid nonparent edges, without endpoint dots
    \foreach \dd in {-.30,-.10,.10}{
      \draw[edge]
        (\name\i)--
        ($(\name\i)+({\dir*\dd},-.85)$);
    }

    % Possible fourth nonparent edge
    \draw[optional]
      (\name\i)--
      ($(\name\i)+({\dir*.30},-.85)$);
  }
}

\end{tikzpicture}

\caption{The rooted configuration at the fixed $5$-vertex $u$.
At each child, the dashed edge is present
if and only if that child has degree five.
The figure is schematic: endpoints at distance three
from $u$ need not be distinct.}
\label{fig:root}
\end{figure}
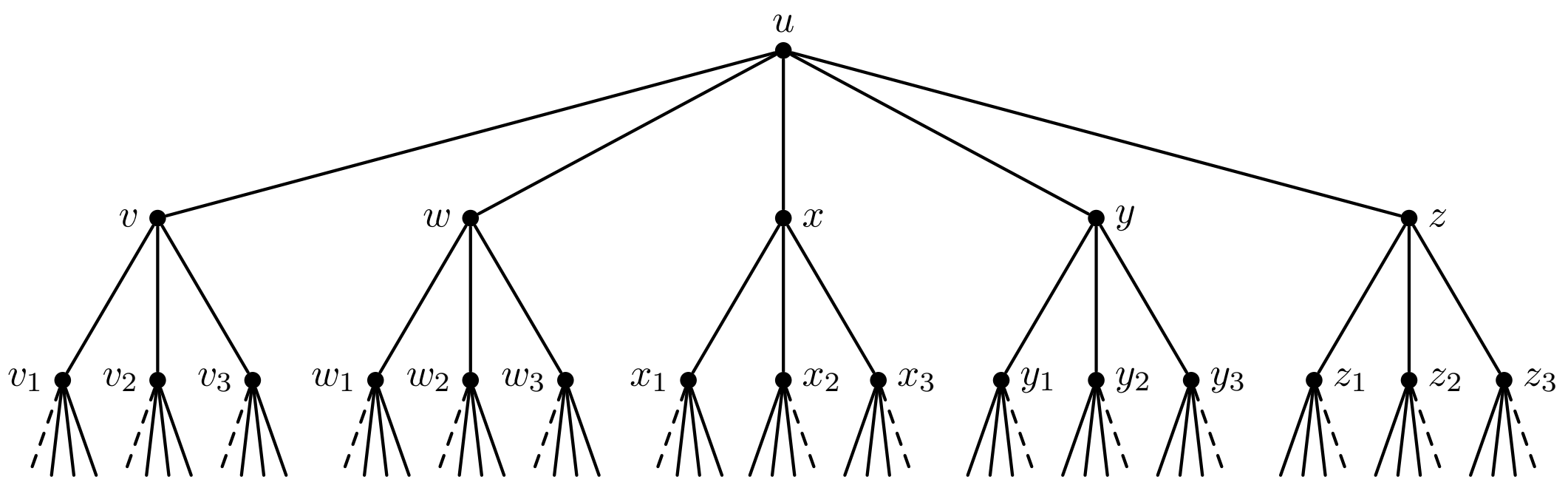

If $G$ had no $5$-vertex, then $\Delta(G)=4$, and the bound
$\chi'_s(G)\le21$ of Huang, Santana, and
Yu~\cite{Huang} would contradict the choice of $G$.
Fix a $5$-vertex $u$ and write $N_G(u)=\{v,w,x,y,z\}$.
Since $\theta(G)\le9$, every neighbor of $u$
has degree four.
We call $v,w,x,y$ the main roots and $uv,uw,ux,uy,uz$
the central edges.
For $a\in\{v,w,x,y,z\}$, write
$N_G(a)=\{u,a_1,a_2,a_3\}$.
We call the vertices $a_i$ the children of $a$ and the edge
$aa_i$ the parent edge of $a_i$.
An edge incident with $a_i$ other than $aa_i$ is called
a nonparent edge at $a_i$.
For each main root $a$, its three children form the branch
rooted at $a$.
For each child $a_i$, write
$N_G(a_i)=\{a,a_{i1},a_{i2},a_{i3},a_{i4}\}$,
where $a_{i4}$ is omitted when $d_G(a_i)=4$.
By \cref{cl:girth}, the vertices
$u,v,w,x,y,z$ and the fifteen children are pairwise distinct.
The rooted configuration is shown in \cref{fig:root}.

\section{An edge ordering and a vertex partition}\label{sec:sequence}

Let $P_0$ be the set of the following ten
edges, precolored as indicated:
\begin{center}
{
\begin{tabular}{cl}
\toprule
Color & Precolored edges\\
\midrule
$1$ & $vv_1,ww_1,xx_1,yy_1$\\
$2$ & $vv_2,ww_2,xx_2$\\
$3$ & $vv_3,ww_3,xx_3$\\
\bottomrule
\end{tabular}
}
\end{center}
Each color class is an induced matching:
a conflict between parent edges in different branches
would create a cycle of length at most five,
contrary to \cref{cl:girth}.

Let $S_0=(y_2y_{21},y_3y_{31},yy_2,yy_3,uv,uw,ux,uz,uy)$.
Choose a longest sequence $S$ of distinct edges of
$E(G)\setminus P_0$ such that $S_0$ is its final segment
and every edge preceding $S_0$ conflicts with at least
four edges occurring later in $S$.
Such a sequence exists, since $S_0$ itself satisfies
these conditions.
In set expressions, a sequence denotes
its set of edges.

\begin{lemma}\label{lem:S}
If $E(G)=P_0\cup S$, then $G$ has a
strong $28$-edge-coloring.
The colors initially assigned to $P_0$ need not all
be preserved.
\end{lemma}

\begin{proof}
Color the edges of $S\setminus S_0$ in their sequence order.
When such an edge $e$ is colored, at least four later
conflicting edges remain uncolored.
Thus \cref{lem:four-savings} applies.

Next try to color the edges of $S_0$ in the stated order.
For the last seven edges, the following
table gives lower bounds $t(e)$ and $r(e)$
as in \cref{lem:certificate}:
\[
\begin{array}{c|ccccccc}
 e&yy_2&yy_3&uv&uw&ux&uz&uy\\ \hline
 t(e)&6&5&4&3&2&1&0\\
 r(e)&0&0&2&2&2&4&4
\end{array}
\]
For $uv,uw,ux$, count the two repetitions
from $vv_2,ww_2,xx_2$; for $uz,uy$, also count the two
from $vv_3,ww_3,xx_3$.
Thus $t(e)+r(e)\ge4$ in every column.

Only $y_2y_{21}$ or $y_3y_{31}$ can therefore
be the first edge with no available color.
Suppose this happens at $e$.
It sees the three uncolored edges
$yy_2,yy_3,uy$.
Since $|\Nstr_G(e)|\le31$, at most $28$ of its strong
neighbors are colored.
Since no color is available, these neighbors use all
$28$ colors exactly once.
In particular, $yy_1$ is its unique strong neighbor
of color $1$.
Uncolor $yy_1$ and assign $1$ to $e$.
If the other of the two edges is still uncolored,
it now sees the four uncolored edges $yy_1,yy_2,yy_3,uy$,
so \cref{lem:four-savings} colors it.
Then color $yy_1$, which sees the seven uncolored edges
$yy_2,yy_3,uv,uw,ux,uz,uy$.
Finish in the original order using the table above,
whose bounds remain valid because the six edges
precolored with $2$ or $3$ retain their colors.
\end{proof}

Put $H=E(G)\setminus(P_0\cup S)$, and let $L$ be the set
of endpoints of edges in $H$.
By \cref{lem:S}, $H\ne\varnothing$.
By the choice of $S$,
\begin{equation}\label{eq:H-three}
 \text{every edge of $H$ conflicts with at most three edges of $S$.}
\end{equation}
Indeed, an edge of $H$ conflicting with four edges of $S$
could be placed at the beginning of the sequence.

None of $u,v,w,x,y,z,y_2,y_3$ belongs to $L$.
For $u,v,w,x,y$, all incident edges belong to $P_0\cup S_0$.
Every edge $zz_i$ conflicts with $uz,uv,uw,ux\in S_0$,
so it belongs to $S$.
Every edge at $y_2$ other than $yy_2$ and $y_2y_{21}$
conflicts with $y_2y_{21},yy_2,yy_3,uy$, and hence belongs
to $S$; the argument for $y_3$ is identical.

We now remove all colors and keep
$P_0,S,H,L$ fixed.
To describe the boundary of $L$, put
$\Azero=\{v_1,v_2,v_3,w_1,w_2,w_3,x_1,x_2,x_3,y_1\}$
and $F=E(L,V(G)\setminus L)$.

\begin{lemma}\label{lem:boundary}
The following statements hold.
\begin{enumerate}[label=(\alph*)]
\item Every edge of $F$ has exactly one endpoint in $\Azero$.
\item Every vertex of $\Azero$ is incident with at most two
      edges of $F$.
\item Every vertex of $L$ is incident with at most one
      edge of $F$.
\item If $p\in\Azero\setminus L$ is incident with
      an edge of $F$, then $d_G(p)=4$.
\end{enumerate}
\end{lemma}

\begin{proof}
Let $e=\ell p\in F$, where $\ell\in L$ and $p\notin L$,
and choose $h\in H$ incident with $\ell$.
Since $p\notin L$, no edge incident with $p$ belongs to $H$.
Thus every edge incident with $p$ outside $P_0$ belongs to $S$.

Suppose that no edge of $P_0$ is incident with $p$.
All of its at least four incident edges then belong to $S$
and conflict with $h$ through $e$, contrary to
\eqref{eq:H-three}.
Therefore $p$ is incident with an edge of $P_0$.
The endpoints of $P_0$ are $v,w,x,y$ and the vertices of
$\Azero$.
If $p\in\{v,w,x,y\}$, then $e$ is a parent edge to a vertex
of $\Azero\cap L$; here we use $y_2,y_3\notin L$.
Otherwise $p\in\Azero$.
Thus $e$ has an endpoint in $\Azero$.
It cannot have two, since an edge between vertices of
$\Azero$ produces a triangle or a $5$-cycle.
This proves (a).

We next prove (c).
Suppose that $\ell\in L$ is incident with at least two
edges of $F$.
If $\ell\notin\Azero$, choose two of them, $\ell p$ and
$\ell q$, and an edge $h\in H$ incident with $\ell$.
By (a), $p,q\in\Azero\setminus L$, and both chosen edges
belong to $S$.
Apart from its parent edge in $P_0$ and $\ell p$, the vertex
$p$ has at least two other incident edges, both in $S$.
These two edges, $\ell p$, and $\ell q$ are distinct and
all conflict with $h$, contrary to \eqref{eq:H-three}.

If instead $\ell\in\Azero$, its parent edge belongs to $F$.
Choose another edge $\ell q\in F$ and an edge $h\in H$
at $\ell$.
By (a), $q\notin\Azero$, and the girth condition excludes
$q\in\{v,w,x,y\}$ because $\ell q$ is not the parent edge.
Hence no edge of $P_0$ is incident with $q$.
All of its at least four incident edges belong to $S$
and conflict with $h$ through $\ell q$, a contradiction.
This proves (c).

For (d), take $p\in\Azero\setminus L$
incident with $p\ell\in F$, and choose $h\in H$ at $\ell$.
The parent edge of $p$ is its only incident edge in $P_0$.
Its other $d_G(p)-1$ incident edges belong to $S$ and
conflict with $h$ through $p\ell$.
Thus \eqref{eq:H-three} gives $d_G(p)-1\le3$.
Since $d_G(p)\ge4$, we have $d_G(p)=4$.

By (c), part (b) holds for vertices of
$\Azero\cap L$.
By (d), every vertex of $\Azero\setminus L$ incident with
$F$ has degree four.
Suppose that such a vertex $p$ has all three nonparent
edges in $F$, and write $f_i=p\ell_i$ for $i\in[3]$,
where $\ell_i\in L$.
Index them so that $f_3$ occurs last in $S$.
Each $f_i$ lies in $S\setminus S_0$.
Choose $h_i\in H$ at $\ell_i$.
Every $h_i$ conflicts with all three edges $f_1,f_2,f_3$.

Let $g$ be a later edge of $S$ conflicting
with $f_3$.
If $g$ is incident with $\ell_1$ or $\ell_2$, then it
conflicts with the corresponding $h_i$.
If $g$ is incident with $\ell_3$ or one of its neighbors,
then it conflicts with $h_3$.
Either possibility contradicts \eqref{eq:H-three}, since
each $h_i$ already conflicts with $f_1,f_2,f_3$.

Since $f_1,f_2$ precede $f_3$, any remaining
$g$ must be incident with the parent of $p$.
If $p\in\Azero\setminus\{y_1\}$, the only edge at its parent
belonging to $S$ is the corresponding central edge.
If $p=y_1$, the only such edges are $yy_2,yy_3,uy$.
Thus $f_3$ has at most three later conflicting edges in $S$,
contrary to the choice of $S$.
This proves (b).
\end{proof}

Define
$V_F^+=\{q\in V(G)\setminus L:q\ell\in F
\text{ for some }\ell\in L\}$,
$M=\{u,v,w,x,y\}\cup V_F^+$, and
$R=V(G)\setminus(L\cup M)$.
By \cref{lem:boundary}(a), an edge of $F$
has an endpoint in $\Azero$.
If this endpoint lies in $L$, its parent edge also belongs
to $F$, so \cref{lem:boundary}(c) forces the two edges
to coincide. Hence
\(
 M\subseteq\{u,v,w,x,y\}\cup\Azero.
\)

\begin{lemma}\label{lem:separation}
We have $E(L,R)=\varnothing$.
Moreover, $z,y_2,y_3,z_1,z_2,z_3\in R$ and
$zz_i\in E(G[R])$ for $i\in[3]$.
\end{lemma}

\begin{proof}
By the definition of $V_F^+$, every neighbor
of $L$ outside $L$ belongs to $M$.
Thus $E(L,R)=\varnothing$.

The containment above and the distinctness
of the rooted vertices imply that none of
$z,y_2,y_3,z_1,z_2,z_3$ belongs to $M$.
We already know from \cref{sec:sequence} that
$z,y_2,y_3\notin L$.
If $z_i\in L$, then $zz_i\in F$ has neither endpoint in
$\Azero$, contrary to \cref{lem:boundary}(a).
Thus all six vertices lie in $R$, and
$zz_i\in E(G[R])$ for $i\in[3]$.
\end{proof}

\begin{lemma}\label{lem:local-rules}
Let $a\in\{v,w,x,y\}$ and $i\in[3]$.
\begin{enumerate}[label=(\alph*)]
\item If $a_i\in L$, then $aa_i\in F$ and every other edge
      at $a_i$ lies in $G[L]$.
\item If $a_i\in R$, then every nonparent edge at $a_i$
      lies in $G[R]$.
\item Every vertex of $M\cap\Azero$ has degree four.
\item If $a_i\in M\cap\Azero$, then
      $\bigl(|E(\{a_i\},L)|,|E(\{a_i\},R)|\bigr)
      \in\{(1,2),(2,1)\}$.
\item If $a_i a_{ij}\in F$, then $a_{ij}\in L$ and every
      other edge at $a_{ij}$ lies in $G[L]$.
\item If $a_i a_{ij}\in E(M,R)$, then $a_{ij}\in R$.
\end{enumerate}
\end{lemma}

\begin{proof}
If $a_i\in L$, its parent lies outside $L$, so $aa_i\in F$.
By \cref{lem:boundary}(c), this is the only edge of $F$
at $a_i$; hence every other edge at $a_i$ lies in $G[L]$.
This proves (a).

If $a_i\in R$, no edge at $a_i$ enters $L$,
by \cref{lem:separation}.
A nonparent edge cannot enter $M$ either: since
$M\subseteq\{u,v,w,x,y\}\cup\Azero$, it would produce a cycle
of length at most five.
Thus (b) holds.

If $a_i\in M\cap\Azero$, then $a_i\in V_F^+$, so it is
incident with an edge of $F$.
Part (c) follows from \cref{lem:boundary}(d).
The parent edge $aa_i$ lies in $G[M]$; each of the other
three edges enters $L$ or $R$, since another edge into $M$
would produce a cycle of length at most five.
At least one and at most two of these edges belong to $F$,
by the definition of $V_F^+$ and \cref{lem:boundary}(b).
This proves (d).

For (e), suppose that $a_i a_{ij}\in F$.
By the girth condition, $a_{ij}\notin\Azero$, so
\cref{lem:boundary}(a) gives
$a_i\in\Azero$.
If $a_i\in L$, its parent edge and $a_i a_{ij}$ would be
two edges of $F$ at $a_i$, contrary to
\cref{lem:boundary}(c).
Hence $a_{ij}\in L$.
The same part of \cref{lem:boundary} implies that every
other edge at $a_{ij}$ lies in $G[L]$.

For (f), part (b) excludes $a_i\in R$.
Since $a_i a_{ij}\in E(M,R)$, we have
$a_i\in M$ and $a_{ij}\in R$.
\end{proof}

\section{Further properties of the partition}
\label{sec:connectors}

We next exclude three configurations.
In each edge ordering below, the edges are colored
from left to right.
An optional fourth nonparent edge is omitted when
the relevant vertex has degree four.
At each coloring step, we apply
\cref{lem:certificate} with the indicated bounds
on $t(e)$ and $r(e)$.

For $q\in R$, let
$\rhoR(q)=d_{G[R]}(q)$.

\begin{lemma}\label{lem:mixed}
No $a\in\{v,w,x,y\}$ has one child in $L$ and another child in $R$.
\end{lemma}

\begin{proof}
Suppose otherwise.
After relabeling, let $a_1\in L$ and $a_3\in R$, and put
$\{b,c,d\}=\{v,w,x,y\}\setminus\{a\}$.
By \cref{lem:local-rules}(a),(b), all nonparent edges at
$a_1$ lie in $G[L]$, whereas those at $a_3$ lie in $G[R]$.

Let $Q_L=G[L]+ma_1$, where $m$ is a new leaf, and let
$Q_R=G[R]+a_3z$.
The edge $a_3z$ does not belong to $G$, since it would create
the $4$-cycle $a_3auza_3$.
In $Q_L$, the degree of $a_1$ equals its degree in $G$,
and the new edge has weight at most $1+5=6$.
In $Q_R$, the degrees of $a_3$ and $z$ equal their degrees
in $G$, and the new edge has weight at most $5+4=9$.
No old edge has weight greater than its weight in $G$.
Both auxiliary graphs are loopless and have
Ore-degree at most nine.
Moreover,
$|V(Q)|+|E(Q)|<|V(G)|+|E(G)|$
for $Q\in\{Q_L,Q_R\}$.

By minimality, color $Q_R$, and permute its colors so that
$zz_1$ has color $4$ and three nonparent edges at $a_3$
have colors $1,2,3$.
Since $zz_1$ sees every edge at $a_3$ through the added
edge $a_3z$, no nonparent edge at $a_3$ has color $4$.
Similarly, color $Q_L$, and permute its colors so that
$ma_1$ has color $4$ and three nonparent edges at $a_1$
have colors $1,2,3$.
The color of $ma_1$ ensures that every edge
of $G[L]$ incident with $a_1$ or one of its neighbors in
$G[L]$ avoids color $4$.
Discard the added edges and retain the colors on $G[L]$
and $G[R]$.
These form a partial strong edge-coloring of $G$, since
both subgraphs are induced and $E(L,R)=\varnothing$.

For each $p\in\{b,c,d\}$, first color every still-uncolored
nonparent edge at $p_1,p_2,p_3$, and then color
$pp_1,pp_2,pp_3$.
Each nonparent edge sees the four reserved edges
$up,pp_1,pp_2,pp_3$, and each parent edge sees all five
still-uncolored central edges.
Thus \cref{lem:four-savings} applies.
During this coloring of the branches rooted at $b,c,d$,
choose a color other than $4$ for each nonparent edge
that conflicts with $aa_1$.
Such an edge sees $aa_1$ in addition to the four reserved
edges, so it has at least two available colors.

If $a_2\in M$, color its still-uncolored nonparent edges
while $ua,aa_1,aa_2,aa_3$ remain uncolored;
\cref{lem:four-savings} applies.
If $a_2\in L$ or $a_2\in R$, its nonparent edges are already
colored.
Only the five central edges and
$aa_1,aa_2,aa_3$ remain uncolored.

\emph{Case 1: color $4$ occurs on a nonparent edge at $a_2$.}
Color
\(
 uz,ub,uc,ud,aa_1,aa_2,aa_3,ua
\)
in this order.
The first four edges see at least $7,6,5,4$ later uncolored
edges, respectively.
Each of $aa_1,aa_2,aa_3$ sees the three repeated colors
$1,2,3$ on the $a_1$- and $a_3$-sides, together with
$3,2,1$ later uncolored edges, respectively.
Finally, $ua$ sees the three repeated colors $1,2,3$ and
also sees color $4$ on both the chosen edge at $a_2$
and $zz_1$.
Thus \cref{lem:four-savings} applies at every step.

\emph{Case 2: no nonparent edge at $a_2$ has color $4$.}
First color $aa_1$ with $4$.
Every colored edge conflicting with $aa_1$
through $a_1$ avoids $4$: retained edges do so by the color
of $ma_1$, and subsequently colored edges do so by the
avoidance rule above.
The colored edges reached through $aa_2$ avoid $4$ by
the case assumption, and those reached through $aa_3$
avoid $4$ by the color of $zz_1$.
All central edges remain uncolored.
Hence $4$ is available for $aa_1$.

Now color
\(
 uz,ub,uc,ud,aa_2,aa_3,ua
\)
in this order.
The first three edges see at least $6,5,4$ later uncolored
edges, respectively.
The edge $ud$ sees the three later edges $ua,aa_2,aa_3$
and sees color $4$ on both $aa_1$ and $zz_1$.
Each of $aa_2,aa_3$ sees the three repeated colors $1,2,3$
and, respectively, two and one later uncolored edges.
The final edge $ua$ sees each of the colors $1,2,3,4$
at least twice.
Again \cref{lem:four-savings} applies at every step.

Both cases yield a strong $28$-edge-coloring
of $G$, a contradiction.
\end{proof}

\begin{lemma}\label{lem:four-branch}
Let $r\in R$ be adjacent to four vertices of $M\cap\Azero$,
one in each of the branches rooted at $v,w,x,y$.
Then $\rhoR(r)\ge1$.
\end{lemma}

\begin{proof}
Suppose that $\rhoR(r)=0$.
Since $E(L,R)=\varnothing$, all neighbors
of $r$ lie in $M$.
The girth condition excludes neighbors in
$\{u,v,w,x,y\}$ and two neighbors in the same branch.
Hence $d_G(r)=4$.
After relabeling the children of $v,w,x$, write
$N_G(r)=\{v_1,w_1,x_1,y_1\}$.
By \cref{lem:local-rules}(c), every neighbor of $r$
has degree four.
Choose an edge $y_1\ell\in F$ with $\ell\in L$,
and write
$N_G(y_1)=\{y,r,\ell,q\}$.
By \cref{lem:boundary}(c), every other edge incident
with $\ell$ lies in $G[L]$.

Let $Q_L=G[L]+m\ell$, where $m$ is a new leaf, and let
$Q_R=G[R]+\{ry_2,ry_3,rz\}$.
The vertices $r,y_2,y_3,z$ are distinct
by \cref{cl:girth}, and the added edges of $Q_R$
do not belong to $G$, since $\rhoR(r)=0$.
The new edge in $Q_L$ has weight at most six.
In $Q_R$, the degree of $r$ is three, whereas the degrees
of $y_2,y_3,z$ are at most their degrees in $G$.
Hence each new edge has weight at most eight, and all
old edge weights remain at most nine.
Both auxiliary graphs are loopless and satisfy
\[
  \theta(Q)\le9,
  \qquad
  |V(Q)|+|E(Q)|<|V(G)|+|E(G)|
  \quad (Q\in\{Q_L,Q_R\}).
\]
By minimality, let $\varphi_R$ be a strong
$28$-edge-coloring of $Q_R$, and permute its colors so that
$\varphi_R(ry_2)=1$, $\varphi_R(ry_3)=2$, and
$\varphi_R(rz)=3$.
Retain only the colors on edges of $G[R]$.
The colors of the added edges allow us to assign colors
$1,2,3$ to $yy_2,yy_3,uz$, respectively:
each assigned color differs from the colors on all
conflicting edges already colored in $G[R]$.
Next color $uy$, and call its color $4$.
At this step, $uy$ sees the four uncolored edges
$uv,uw,ux,yy_1$, so \cref{lem:four-savings} applies.
Since $uy$ also sees $yy_2,yy_3,uz$,
its color differs from $1,2,3$.

By minimality, color $Q_L$.
The added edge $m\ell$ and any three edges of $G[L]$
incident with $\ell$ are pairwise conflicting.
Permute the colors so that $m\ell$ has color $3$ and
the three chosen edges have colors $1,2,4$.
Retain only the colors on edges of $G[L]$.
The color of $m\ell$ ensures that every
edge of $G[L]$ incident with $\ell$ or one of its
neighbors in $G[L]$ avoids color $3$.
The retained colors, together with the four assigned colors,
form a partial strong edge-coloring of $G$:
$E(L,R)=\varnothing$, and none of $yy_2,yy_3,uz,uy$
conflicts with an edge of $G[L]$.

In all greedy steps below preceding $v_1r$,
avoid color $3$ on any edge conflicting with $y_1\ell$.
The bounds below give each such edge at least four other
uncolored strong neighbors, so it has at least two
available colors.

For each $p\in\{v,w,x\}$, first color every still-uncolored
nonparent edge at $p_2$ and $p_3$.
Each such edge sees the four reserved edges
$up,pp_1,pp_2,pp_3$.
Put $X_p=N_G(p_1)\setminus\{p,r\}$; then $|X_p|=2$.
Color the remaining edges in the order
\[
\begin{aligned}
&ux,xx_3,xx_2,xx_1,\{x_1s:s\in X_x\},\\
&uw,ww_3,ww_2,ww_1,\{w_1s:s\in X_w\},\\
&uv,vv_3,vv_2,vv_1,\{v_1s:s\in X_v\},\\
&v_1r,w_1r,x_1r,y_1\ell,y_1q,y_1r,yy_1.
\end{aligned}
\]
The edges within each set of braces may be colored
in any order.

Every edge before $v_1r$ in this ordering has at least
four later uncolored conflicting edges.
Indeed, $up$ sees $pp_1,pp_2,pp_3,p_1r$;
each of $pp_3,pp_2$ sees $pp_1,p_1r$ and the two edges
$p_1s$ with $s\in X_p$; and both $pp_1$ and each edge
$p_1s$ see all four edges incident with $r$.
The edges $v_1r,w_1r,x_1r$ have at least six, five, and
four later uncolored conflicting edges, respectively.
Thus \cref{lem:four-savings} applies at every step
through $x_1r$.

After coloring $x_1r$, consider whether color $3$ appears
on one of $v_1r,w_1r,x_1r$ or on a colored edge incident
with $q$.
If so, color $y_1\ell$ using \cref{lem:four-savings}.
Indeed, it sees the three later uncolored edges
$y_1q,y_1r,yy_1$.
In addition, the colors $1,2,4$ each occur at least twice
among $yy_2,yy_3,uy$ and the three selected edges
incident with $\ell$.

If color $3$ appears in none of these positions,
assign color $3$ to $y_1\ell$.
This assignment is valid: colored strong neighbors
reached through $\ell$ avoid $3$ by the auxiliary coloring
and the avoidance rule; colored edges incident with
$r$ or $q$ avoid $3$ by assumption; and the colored
edges at $y$ have colors $1,2,4$.
All remaining edges incident with $y_1$ are uncolored.

Next color $y_1q$ and then $y_1r$.
These edges see, respectively, two and one later uncolored
edges, together with the three repeated colors $1,2,4$.
Thus \cref{lem:four-savings} applies to both edges.

It remains to color $yy_1$.
The colors $1,2,4$ each occur at least twice in its
colored strong neighborhood.
Color $3$ also occurs at least twice: once on $uz$
and once on one of $v_1r,w_1r,x_1r$, on an edge
incident with $q$, or on $y_1\ell$.
Hence $r(yy_1)\ge4$, and \cref{lem:four-savings}
gives an available color for $yy_1$.
This completes a strong $28$-edge-coloring of $G$,
a contradiction.
\end{proof}

\begin{lemma}\label{lem:M-nonempty}
We have $M\cap\Azero\ne\varnothing$.
\end{lemma}

\begin{proof}
Suppose that $M\cap\Azero=\varnothing$.
Since $y_2,y_3\in R$, \cref{lem:mixed} implies that
$y_1\notin L$, and hence all three children of $y$ lie in $R$.
For each $a\in\{v,w,x\}$, either all three children lie in $L$
or all three lie in $R$.

Since $H\ne\varnothing$, the set $L$ is nonempty.
Since $u\notin L$ and $G$ is connected,
we have $F\ne\varnothing$.
By \cref{lem:boundary}(a), every edge of $F$ has an endpoint
in $\Azero$.
This endpoint lies in $L$, since otherwise it would belong
to $M\cap\Azero$.
Its parent edge also belongs to $F$, so
\cref{lem:boundary}(c) forces the boundary edge to be
its parent edge.

Let $B\subseteq\{v,w,x\}$ be the set of roots whose three
children lie in $L$.
Then $B\ne\varnothing$, and
\[
  E(L\cup B,V(G)\setminus(L\cup B))
    =\{ua:a\in B\}.
\]
This is an edge-cut of size at most three,
contrary to \cref{cl:cut}.
\end{proof}

The following lemma gives the preliminary coloring steps
used in the four cases of \cref{lem:low-return}.

\begin{lemma}\label{lem:preliminary-sweep}
Suppose that a partial strong
$28$-edge-coloring of $G$ has been obtained by coloring
$G[L]$ and $G[R]$ and making the additional assignments
specified in one of Cases~1--4 in the proof of
\cref{lem:low-return}.
All remaining uncolored edges, except the central edges
and those included in the final ordering of that case,
can be colored in the following two stages.
\begin{enumerate}[label=(\roman*)]
\item For each $a_i\in M\cap\Azero$, first color every
still-uncolored nonparent edge at $a_i$ that is not
included in the final ordering, while
$ua,aa_1,aa_2,aa_3$ remain uncolored.
\item Then color every still-uncolored parent edge $aa_i$
not included in the final ordering, while at least four
of the five central edges $uv,uw,ux,uy,uz$ remain uncolored.
\end{enumerate}
Every step extends the current partial strong
$28$-edge-coloring.
\end{lemma}

\begin{proof}
Every edge colored in (i) conflicts with
$ua,aa_1,aa_2,aa_3$, and every parent edge colored in (ii)
conflicts with all five central edges.
Thus at least four strong neighbors remain uncolored
at each step, and \cref{lem:four-savings} applies.
\end{proof}

The following bounds will be used for the last edges
in the four coloring orders.

\begin{lemma}\label{lem:low-return-bounds}
Let $p\in M\cap\Azero$, let $a$ be its parent,
and let $q\in N_G(p)\cap R$.
\begin{enumerate}[label=(\alph*)]
\item If $d_G(q)=4$, $\rhoR(q)=1$, and the two neighbors
of $q$ other than $p$ and its unique $R$-neighbor are
$4$-vertices $c,d\in M\cap\Azero$, then
$|\Nstr_G(pq)|\le27$ and $|\Nstr_G(cq)|\le27$.
\item If $\rhoR(q)=2$, then $|\Nstr_G(pq)|\le30$.
More precisely, the bound is at most $28$ when $d_G(q)=4$
and at most $30$ when $d_G(q)=5$.
\item We have $|\Nstr_G(ap)|\le30$,
with the stronger bound $29$ when $d_G(q)=4$.
\end{enumerate}
\end{lemma}

\begin{proof}
For any edge $xy$, count first the edges incident with
$x$ or $y$ and then the edges reached through their
other neighbors.
This gives
\begin{equation}\label{eq:refined-local-count}
 |\Nstr_G(xy)|\le d_G(x)+d_G(y)-2
 +\sum_{z\in N_G(x)\setminus\{y\}}(d_G(z)-1)
 +\sum_{z\in N_G(y)\setminus\{x\}}(d_G(z)-1).
\end{equation}
Possible repetitions only decrease the number of
distinct strong neighbors.

By \cref{lem:local-rules}(c), $d_G(p)=4$,
and its parent $a$ also has degree four.
Thus, for $pq$, the three neighbors of $p$ other than $q$
contribute at most $3+4+4=11$.

For (a), the three neighbors of $q$
other than $p$ contribute at most $3+3+4=10$.
Hence \eqref{eq:refined-local-count} gives
$|\Nstr_G(pq)|\le6+11+10=27$.
The same argument applies to $cq$, since
$c\in M\cap\Azero$.

For (b), suppose first that $d_G(q)=4$.
Since $q$ is adjacent to $p$, the girth
condition excludes neighbors of $q$ in $\{u,v,w,x,y\}$.
Together with $E(L,R)=\varnothing$ and the description
of $M$, this implies that every neighbor of $q$ outside
$R$ belongs to $M\cap\Azero$.
Besides $p$, the vertex $q$ therefore has one further
neighbor in $M\cap\Azero$ and two neighbors in $R$.
Their contributions are at most $3+4+4$, so
\eqref{eq:refined-local-count} gives
$|\Nstr_G(pq)|\le6+11+11=28$.
If $d_G(q)=5$, every neighbor of $q$ has degree at most
four, and hence
$|\Nstr_G(pq)|\le7+11+4\cdot3=30$.

Finally, since $ap$ is a $(4,4)$-edge, the bound
$|\Nstr_G(ap)|\le30$ follows from
\cref{lem:strong-neighborhood}.
If $d_G(q)=4$, then for $ap$, the three neighbors of $p$
other than $a$ contribute at most $3+4+4=11$,
and the three neighbors of $a$ other than $p$
contribute at most $12$.
Thus $|\Nstr_G(ap)|\le6+11+12=29$.
\end{proof}

The following lemma shows how to transfer
colors from added edges to edges of $G$.
Initially, only the edges of $G[L]$ and $G[R]$ are colored,
and all edges incident with $M$ remain uncolored.

\begin{lemma}\label{lem:witness-transfer}
In each part below, suppose that $G[L]$
and the indicated auxiliary graph have strong
$28$-edge-colorings.
Retain only their colors on $G[L]$ and $G[R]$,
leaving all edges incident with $M$ uncolored.
\begin{enumerate}[label=(\alph*)]
\item Let $p=y_1\in M$, let $q\in R$ be adjacent to $p$,
and add the two edges $qy_2,qy_3$ to $G[R]$.
If $qy_i$ has color $c_i$ in the auxiliary coloring,
then $yy_i$ may be assigned color $c_i$, for $i=2,3$.
\item Let $p=v_1\in M$, let $q\in R$ be adjacent to $p$,
and add the edge $qz$ to $G[R]$.
If $qz$ has color $c$ in the auxiliary coloring,
then $uz$ may be assigned color $c$.
\end{enumerate}
In either part, these assignments extend
the retained coloring to a partial strong
$28$-edge-coloring of $G$.
\end{lemma}

\begin{proof}
The retained colors form a partial strong
edge-coloring of $G$, since $G[L]$ and $G[R]$ are induced
and $E(L,R)=\varnothing$.

For (a), consider a colored edge of $G[R]$ that conflicts
with $yy_2$.
Such an edge is incident with $y_2$, with a neighbor
of $y_2$ in $R$, or with $y_3$.
The added edge $qy_2$ sees the first two
types through $y_2$ and the third through $qy_3$.
Thus the color of $qy_2$ is absent from every colored edge
of $G[R]$ that conflicts with $yy_2$.
The argument for $yy_3$ is symmetric.
No edge of $G[L]$ conflicts with $yy_i$,
since $N_G(y)\cap L=\varnothing$
and $E(L,R)=\varnothing$.
The added edges $qy_2,qy_3$ are adjacent, so $c_2\ne c_3$.
Thus the assignments to $yy_2,yy_3$ are compatible
with each other and with the already colored edges.

For (b), every colored edge of $G[R]$ that conflicts
with $uz$ is incident with $z$ or with a neighbor
of $z$ in $R$.
The added edge $qz$ sees every such edge.
No colored edge of $G[L]$ conflicts with $uz$,
since $N_G(u)\cap L=\varnothing$
and $E(L,R)=\varnothing$.
Hence the color of $qz$ is available for $uz$.
\end{proof}

\begin{lemma}\label{lem:low-return}
Assume that
\begin{equation}\label{eq:all-low-return}
   \sum_{q\in N_G(p)\cap R}\rhoR(q)\le2
   \qquad\text{for every }p\in M\cap\Azero.
\end{equation}
Then $G$ admits a strong $28$-edge-coloring.
Consequently, \eqref{eq:all-low-return} cannot hold in a
minimum counterexample.
\end{lemma}

\begin{proof}
By \cref{lem:local-rules}(c), every vertex in
$M\cap\Azero$ has degree four.

Let $q\in R$ be adjacent to a vertex of
$M\cap\Azero$.
The girth condition excludes neighbors of $q$ in
$\{u,v,w,x,y\}$ and permits at most one neighbor in
$\Azero$ from each main branch.
Since $E(L,R)=\varnothing$, every neighbor of $q$ outside
$R$ therefore belongs to $M\cap\Azero$.
If $\rhoR(q)=0$, then $d_G(q)\ge4$ forces $q$ to have a
neighbor in each of the four main branches, contrary to
\cref{lem:four-branch}.
Thus $\rhoR(q)\ge1$, and \eqref{eq:all-low-return}
gives $\rhoR(q)\in\{1,2\}$.

Whenever $p\in M\cap\Azero$ is selected, choose an edge
$p\ell\in F$ with $\ell\in L$.
By \cref{lem:boundary}(c), this is the only boundary edge
at $\ell$.
Since $d_G(\ell)\ge4$, at least three edges of $G[L]$
are incident with $\ell$; denote three of them by
$g_1,g_2,g_3$.
These edges are pairwise adjacent and therefore receive
distinct colors in every strong coloring of $G[L]$.
By permuting the colors, we may assign any prescribed
three distinct colors to these edges.

In each case, color $G[L]$ and the indicated auxiliary
graph $Q_R$ obtained from $G[R]$ by minimality.
We use $\varphi$ for the resulting
colorings and their extensions.
No edge of $G[L]$ conflicts with an edge of $G[R]$:
a line-graph path of length two between them would have
an intermediate edge in $E(L,R)$, which is empty.
The retained colorings can therefore be combined.
Each auxiliary graph is loopless and has
Ore-degree at most nine, as checked below.
It omits at least six vertices of $G$, namely
$u,v,w,x,y$ and a vertex of $L$, and adds at most two edges.
Hence
$|V(Q_R)|+|E(Q_R)|<|V(G)|+|E(G)|$.
The degrees in $Q_R$ do not exceed the corresponding
degrees in $G$, so old edge weights do not increase.
We use colors of added edges only as specified in
\cref{lem:witness-transfer}.

In each case, first choose and permute
the colorings and transfer the specified colors.
Then apply \cref{lem:preliminary-sweep}, preserving
the edges reserved for the final ordering.
All relabelings below involve only $v,w,x$ and their
children.

\smallskip
\noindent\textbf{Case 1: $y_1\in M$ and a chosen
$R$-neighbor $q$ of $y_1$ has $\rhoR(q)=1$.}
Put $p=y_1$, let $f_1$ be the unique edge of $G[R]$
incident with $q$, and write
$N_G(p)\setminus\{y\}=\{\ell,s,q\}$.

Besides $p$ and its unique $R$-neighbor, $q$ has two or
three neighbors in $M\cap\Azero$, lying in distinct
branches among those rooted at $v,w,x$.
Relabel so that $w_1q,x_1q\in E(G)$.
Since every vertex of $M\cap\Azero$ has degree four,
write
$N_G(w_1)=\{w,q,r_1,r_2\}$ and
$N_G(x_1)=\{x,q,s_1,s_2\}$.

The girth condition gives
$q\notin\{y_2,y_3\}$ and $qy_2,qy_3\notin E(G)$,
since either added edge would create a $4$-cycle
through $q,p,y$.
Let $Q_R=G[R]+\{qy_2,qy_3\}$.
Here $d_{Q_R}(q)=3$, and the degrees of $y_2,y_3$
are at most their degrees in $G$.
Thus each new edge has weight at most eight.
Color $Q_R$, and permute its colors so that
$\varphi(f_1)=1$, $\varphi(qy_2)=2$, and
$\varphi(qy_3)=3$.
By \cref{lem:witness-transfer}, assign
$\varphi(yy_2)=2$ and $\varphi(yy_3)=3$.
Permute the coloring of $G[L]$ so that
$\varphi(g_i)=i$ for $i=1,2,3$.
The members of each of the same-color
pairs $(g_1,f_1)$, $(g_2,yy_2)$, and $(g_3,yy_3)$
are nonconflicting by $E(L,R)=\varnothing$ and
\cref{lem:witness-transfer}(a).

Apply \cref{lem:preliminary-sweep} to all
remaining uncolored edges except the central edges and
those listed below, including any additional edge from
$q$ to $M\cap\Azero$.
Then color the listed edges in order:
\begin{equation}\label{eq:low-return-order-1}
\begin{aligned}
&uz,uv,uw,ux,
 ww_3,ww_2,w_1r_1,w_1r_2,ww_1,\\
&xx_3,xx_2,x_1s_1,x_1s_2,xx_1,
 uy,w_1q,x_1q,p\ell,ps,pq,yp.
\end{aligned}
\end{equation}
The following array gives four later uncolored strong
neighbors for every edge preceding $p\ell$:
\[
\begin{array}{c|l@{\qquad}c|l}
uz&uv,uw,ux,uy
&uv&uw,ux,uy,ww_3\\
uw&ux,uy,ww_3,ww_2
&ux&uy,xx_3,xx_2,xx_1\\
ww_3&ww_2,ww_1,w_1r_1,w_1r_2
&ww_2&ww_1,w_1r_1,w_1r_2,w_1q\\
w_1r_1&w_1r_2,ww_1,w_1q,x_1q
&w_1r_2&ww_1,w_1q,x_1q,pq\\
ww_1&uy,w_1q,x_1q,pq
&xx_3&xx_2,xx_1,x_1s_1,x_1s_2\\
xx_2&xx_1,x_1s_1,x_1s_2,x_1q
&x_1s_1&x_1s_2,xx_1,x_1q,w_1q\\
x_1s_2&xx_1,x_1q,w_1q,pq
&xx_1&uy,x_1q,w_1q,pq\\
uy&p\ell,ps,pq,yp
&w_1q&x_1q,pq,p\ell,ps\\
x_1q&pq,p\ell,ps,yp&&
\end{array}
\]
Hence all these edges are colorable by
\cref{lem:four-savings}.
Each of $p\ell,ps,pq,yp$ sees the three same-color pairs
listed above.
Their numbers of later uncolored conflicting edges are
$3,2,1,0$, respectively, and $yp$ is a $(4,4)$-edge.
The edges are therefore colorable in the stated order by
\cref{lem:four-savings,lem:three-savings}.

\smallskip
\noindent\textbf{Case 2: $y_1\in M$ and a chosen
$R$-neighbor $q$ of $y_1$ has $\rhoR(q)=2$.}
Again put $p=y_1$ and
$N_G(p)\setminus\{y\}=\{\ell,s,q\}$.
Let $f_1,f_2$ be the two edges of $G[R]$ incident with $q$.
The vertex $q$ has at least one additional neighbor in
$M\cap\Azero$; after relabeling, let this neighbor be $v_1$.
Write $N_G(v_1)=\{v,q,r_1,r_2\}$.
The edges $qy_2,qy_3$ do not belong to $G$, by the same
$4$-cycle argument as in Case~1.
Let $Q_R=G[R]+\{qy_2,qy_3\}$.
Here $d_{Q_R}(q)=4$, so each new edge has weight at most nine.
Color $Q_R$, and permute its colors so that
$\varphi(f_1)=1$, $\varphi(f_2)=2$, and
$\varphi(qy_2)=3$.
Assign $\varphi(yy_2)=3$ by \cref{lem:witness-transfer},
and permute the colors on $G[L]$ so that
$\varphi(g_i)=i$.
The three same-color pairs are
$(g_1,f_1)$, $(g_2,f_2)$, and $(g_3,yy_2)$.
Their members are nonconflicting for the same reasons
as in Case~1.

Apply \cref{lem:preliminary-sweep} to all
remaining uncolored edges except the central edges and
those listed below, including any additional edge from
$q$ to $M\cap\Azero$.
Then color
\begin{equation}\label{eq:low-return-order-2}
\begin{aligned}
&v_1r_1,v_1r_2,vv_3,vv_2,vv_1,
 uz,ux,uw,uv,uy,\\
&v_1q,yy_3,p\ell,ps,pq,yp.
\end{aligned}
\end{equation}
For each edge preceding $p\ell$, the following table
lists four later uncolored edges that conflict with it:
\[
\begin{array}{c|l@{\qquad}c|l}
v_1r_1&v_1r_2,vv_1,v_1q,pq
&v_1r_2&vv_1,v_1q,pq,vv_2\\
vv_3&vv_2,vv_1,v_1q,uz
&vv_2&vv_1,v_1q,uz,ux\\
vv_1&v_1q,uz,ux,uw
&uz&ux,uw,uv,uy\\
ux&uw,uv,uy,yy_3
&uw&uv,uy,yy_3,yp\\
uv&uy,v_1q,yy_3,yp
&uy&yy_3,p\ell,ps,pq\\
v_1q&pq,p\ell,ps,yp
&yy_3&yp,p\ell,ps,pq
\end{array}
\]
Each of the last four edges $p\ell,ps,pq,yp$ sees the
three same-color pairs above, and their numbers of later
uncolored conflicting edges are $3,2,1,0$, respectively.
Since $yp$ is a $(4,4)$-edge, the extension follows from
\cref{lem:four-savings,lem:three-savings}.

\smallskip
\noindent\textbf{Case 3: $y_1\notin M$, and some vertex
$q\in R$ adjacent to $M\cap\Azero$ has $\rhoR(q)=2$.}
By \cref{lem:mixed} and the fact that $y_2,y_3\in R$,
we have $y_1\notin L$ and hence $y_1\in R$.
Choose an edge $pq$ with $p\in M\cap\Azero$, $q\in R$,
and $\rhoR(q)=2$.
The vertex $p$ lies in a branch rooted at one of $v,w,x$;
after relabeling, write $p=v_1$.
Put $N_G(p)\setminus\{v\}=\{\ell,s,q\}$.
Let $f_1,f_2$ be the two edges of $G[R]$ incident with $q$.

Let $Q_R=G[R]+qz$.
The girth condition gives $q\ne z$,
and $qz\notin E(G)$ since it would create the
$5$-cycle $qpvuzq$.
In $Q_R$, we have $d_{Q_R}(q)=3$ and $d_{Q_R}(z)=4$,
so the new edge has weight seven.
Color $Q_R$, and permute its colors so that
$\varphi(f_1)=1$, $\varphi(f_2)=2$, and $\varphi(qz)=3$.
Assign $\varphi(uz)=3$ by \cref{lem:witness-transfer},
and permute the colors on $G[L]$ so that
$\varphi(g_i)=i$.
The members of each of the same-color
pairs $(g_1,f_1)$, $(g_2,f_2)$, and $(g_3,uz)$
are nonconflicting by $E(L,R)=\varnothing$ and
\cref{lem:witness-transfer}(b).

Apply \cref{lem:preliminary-sweep} to all
remaining uncolored edges except the central edges and
those listed below, including every edge from $q$ to
$M\cap\Azero$ other than $pq$.
Then use the order
\begin{equation}\label{eq:low-return-order-3}
   uy,ux,uw,uv,vv_3,vv_2,p\ell,ps,pq,vp.
\end{equation}
The following array lists four later uncolored
conflicting edges for each of the first six edges:
\[
\begin{array}{c|l@{\qquad}c|l}
uy&ux,uw,uv,vp
&ux&uw,uv,vp,vv_3\\
uw&uv,vp,vv_3,vv_2
&uv&vv_3,vv_2,vp,p\ell\\
vv_3&vv_2,vp,p\ell,ps
&vv_2&vp,p\ell,ps,pq
\end{array}
\]
Each of $p\ell,ps,pq$ sees the two repeated colors on
$(g_1,f_1)$ and $(g_2,f_2)$ and has, respectively,
$3,2,1$ later uncolored conflicting edges.
For $pq$, in particular,
\cref{lem:low-return-bounds}(b) shows that at most
$30-1-2=27$ distinct colors are forbidden.
The final edge $vp$ sees these two repeated colors
and also color $3$ on $g_3$ and $uz$.
Since $vp$ is a $(4,4)$-edge, it is colorable by
\cref{lem:three-savings}.

\smallskip
\noindent\textbf{Case 4: $y_1\notin M$, and every vertex
of $R$ adjacent to $M\cap\Azero$ has degree one in $G[R]$.}
As in Case~3, \cref{lem:mixed} gives $y_1\in R$.
Choose an edge $pq$ with $p\in M\cap\Azero$, $q\in R$,
and $\rhoR(q)=1$.
After relabeling, let $p=v_1$, and write
$N_G(p)\setminus\{v\}=\{\ell,s,q\}$.
The vertex $q$ cannot have degree five.
Indeed, a $5$-vertex with degree one in $G[R]$ would
need three further neighbors in $M\cap\Azero$ besides $p$.
The girth condition permits at most one such neighbor
in each main branch, whereas all three children of $y$
lie in $R$; only the branches rooted at $w$ and $x$ remain.
Hence $d_G(q)=4$, and its two additional neighbors are
$4$-vertices in two distinct branches, which we denote
by $w_1$ and $x_1$.
Let $f_1$ be the unique edge of $G[R]$ incident with $q$.

Let $Q_R=G[R]+qz$.
As in Case~3, $q\ne z$ and
$qz\notin E(G)$.
Here $d_{Q_R}(q)=2$ and $d_{Q_R}(z)=4$,
so the new edge has weight six.
By minimality, color $Q_R$ and permute its colors so that
$\varphi(f_1)=1$ and $\varphi(qz)=3$.
By \cref{lem:witness-transfer}, assign $\varphi(uz)=3$.
Permute the colors on $G[L]$ so that $\varphi(g_i)=i$.
As in Case~3, each of $(g_1,f_1)$ and $(g_3,uz)$ consists
of two nonconflicting edges of the same color.

Apply \cref{lem:preliminary-sweep} to all
remaining uncolored edges except the central edges and
those listed below.
In particular, color $x_1q$ and leave $w_1q$ uncolored.
Then use the order
\begin{equation}\label{eq:low-return-order-4}
  uy,ux,uw,uv,vv_3,vv_2,p\ell,ps,w_1q,vp,pq.
\end{equation}
For the first six edges, the four later uncolored
conflicting edges are the same as in Case~3.
The edge $p\ell$ sees the four later edges
$ps,w_1q,vp,pq$.
The edge $ps$ sees the three later edges $w_1q,vp,pq$
and the repeated color on $(g_1,f_1)$.
By \cref{lem:low-return-bounds}(a), the edge $w_1q$
has at most $27$ strong neighbors and is therefore
colorable when reached.
The edge $vp$ sees the still-uncolored edge $pq$
and two repeated colors, one on $(g_1,f_1)$ and
the other on $(g_3,uz)$.
By \cref{lem:low-return-bounds}(c),
at most $29-1-2=26$ distinct colors are forbidden
for $vp$.
Finally, $pq$ is colorable by its bound $27$.

These cases are exhaustive.
If $y_1\in M$, then \cref{lem:local-rules}(d) gives an
$R$-neighbor of $y_1$, whose degree in $G[R]$ is one
or two; hence Case~1 or Case~2 applies.
If $y_1\notin M$, then $y_1\in R$, as shown above.
By \cref{lem:M-nonempty}, some vertex of $M\cap\Azero$
lies in a branch rooted at $v,w,x$;
by \cref{lem:local-rules}(d), it has an $R$-neighbor.
If any such neighbor has degree two in $G[R]$,
Case~3 applies; otherwise Case~4 applies.
Each case yields a strong $28$-edge-coloring of $G$,
a contradiction.
\end{proof}

\begin{lemma}\label{lem:terminal-vertex}
There exists $p\in M\cap\Azero$ such that
\begin{equation}\label{eq:terminal-good}
   \sum_{q\in N_G(p)\cap R}\rhoR(q)\ge3.
\end{equation}
\end{lemma}

\begin{proof}
If no such vertex existed, then
\eqref{eq:all-low-return} would hold, since each sum
is an integer.
By \cref{lem:low-return}, $G$ would have a strong
$28$-edge-coloring, a contradiction.
\end{proof}

\section{Completing the coloring}

Choose $p=a_1\in M\cap\Azero$ satisfying
\eqref{eq:terminal-good}, where $a\in\{v,w,x,y\}$,
relabeling within its branch if necessary.
By \cref{lem:local-rules}(c), $d_G(p)=4$,
and part (d) gives a neighbor $q_0\in L$
joined to $p$ by a nonparent edge; fix one such vertex.

By \eqref{eq:terminal-good}, at least three distinct edges
of $G[R]$ are incident with vertices of $N_G(p)\cap R$.
Indeed, no edge of $G[R]$ joins two vertices of
$N_G(p)\cap R$, since such an edge would form a triangle
with $p$.
Choose distinct edges $f_1,f_2,f_3\in E(G[R])$,
each incident with a vertex of $N_G(p)\cap R$.
The edge $pq_0$ is the only boundary edge at $q_0$
by \cref{lem:boundary}(c).
Since $d_G(q_0)\ge4$, at least three edges of $G[L]$
are incident with $q_0$; fix three of them and denote
them by $g_1,g_2,g_3$.

\begin{lemma}\label{lem:synchronization}
The graphs $G[L]$ and $G[R]$ have strong
$28$-edge-colorings whose union is a partial strong
edge-coloring $\varphi$ of $G$ satisfying
$\varphi(g_i)=i$ for $1\le i\le3$ and
$\{\varphi(f_1),\varphi(f_2),\varphi(f_3)\}
\subseteq\{1,2,3\}$.
\end{lemma}

\begin{proof}
By minimality, both $G[L]$ and $G[R]$ have strong
$28$-edge-colorings.
The edges $g_1,g_2,g_3$ are incident with $q_0$,
so their colors are distinct.
Permute the colors on $G[L]$ so that they receive
colors $1,2,3$.
The edges $f_1,f_2,f_3$ use at most three colors;
permute the colors on $G[R]$ so that their colors
belong to $\{1,2,3\}$.

Each coloring is a partial strong
edge-coloring of $G$, since its subgraph is induced.
No edge of $G[L]$ conflicts with an edge of $G[R]$:
a path of length two between them in $L(G)$ would have
an intermediate edge in $E(L,R)$, which is empty
by \cref{lem:separation}.
Hence the two colorings can be combined.
\end{proof}

\begin{lemma}\label{lem:terminal-extension}
The partial coloring given by \cref{lem:synchronization}
extends to a strong $28$-edge-coloring of $G$.
\end{lemma}

\begin{proof}
Retain the colorings of $G[L]$ and $G[R]$ given by
\cref{lem:synchronization}.
Every remaining uncolored edge has at least one
endpoint in $M$.

First color every uncolored nonparent edge at a vertex
of $(M\cap\Azero)\setminus\{p\}$.
If this vertex is $b_i$, where $b\in\{v,w,x,y\}$,
the edge conflicts with the four still-uncolored edges
$ub,bb_1,bb_2,bb_3$.
Thus \cref{lem:four-savings} applies.
Leave the three nonparent edges at $p$ uncolored.

Next color every still-uncolored parent edge other than
$aa_1,aa_2,aa_3$.
Each conflicts with all five central edges, which remain
uncolored, so \cref{lem:four-savings} applies.
Then color every central edge other than $ua$.
Each conflicts with the four still-uncolored edges
$ua,aa_1,aa_2,aa_3$, so the same lemma applies.

Now color $ua,aa_2,aa_3$ in this order.
The three nonparent edges at $p=a_1$ are still uncolored.
Together with $aa_1$ and the later members of this order,
they give at least six, five, and four later uncolored
conflicting edges, respectively.
Thus \cref{lem:four-savings} applies
to each edge.

Only the three nonparent edges at $p$ and $aa_1=ap$
remain uncolored.
All six edges $g_1,g_2,g_3,f_1,f_2,f_3$
are incident with neighbors of $p$, so every remaining
edge sees all of them.
These six distinct edges use exactly three colors,
and hence $r_\varphi(e)\ge6-3=3$ for every remaining
edge $e$.
Their colors remain unchanged throughout the extension.

Color the three nonparent edges at $p$ in any order.
Each still conflicts with the uncolored
edge $ap$, so we may take $t(e)=1$ and $r(e)=3$.
Thus \cref{lem:four-savings} applies.

Finally, $ap$ is a $(4,4)$-edge and still has
$r_\varphi(ap)\ge3$.
By \cref{lem:three-savings},
an available color exists for $ap$.
This completes the coloring of $G$.
\end{proof}

\begin{proof}[Proof of Theorem~\ref{thm:main}]

Let $G$ be the minimum counterexample fixed in
\cref{sec:minimal}.
By \cref{lem:terminal-vertex}, there exists
$p\in M\cap\Azero$ satisfying \eqref{eq:terminal-good}.
The partial strong edge-coloring supplied by
\cref{lem:synchronization} extends to a strong
$28$-edge-coloring of $G$ by \cref{lem:terminal-extension},
contrary to the choice of $G$.
This proves \cref{thm:main}.
\end{proof}

\appendix

\section{Short-cycle reductions}\label{app:short-cycles}

We prove the statements used in \cref{cl:girth}
by deleting vertices, coloring the remaining induced
subgraph by minimality, and extending the resulting
partial strong edge-coloring of $G$.

For structural properties of $G$,
we use only \cref{cl:simple,cl:degrees,cl:cut} and
conclusions proved earlier in this appendix;
in particular, we do not assume \cref{cl:girth}.
All conflicts and available colors are measured in
the original graph $G$.
Whenever a coloring is transferred to $G$, we verify
that the retained colors form a partial strong
edge-coloring of $G$.

\subsection{Triangles and complete bipartite subgraphs}

\begin{lemma}\label{lem:no-triangle-app}
The graph $G$ contains no triangle.
\end{lemma}

\begin{proof}

A triangle contains at most one $5$-vertex.
Suppose that $uvw$ is a triangle with $d_G(u)=4$,
and write $N_G(u)=\{v,w,x,y\}$.
By minimality, color $G-u$, obtaining a partial
strong edge-coloring of $G$.

For $ux$, counting through $x$ and then at $v,w,y$
gives at most
\[
  (d_G(x)-1)(9-d_G(x))+3+4+4-1\le26
\]
colored conflicting edges, where the subtraction
accounts for $vw$ counted twice.
Thus $|A(ux)|\ge2$, and similarly $|A(uy)|\ge2$.

For $uv$, the retained edges at $v$ and the external
edges at $w$ contribute at most
$d_G(v)+d_G(w)-3\le6$.
The edges at $x,y$ contribute at most eight, and
those reached through the external neighbors of $v$
contribute at most
$(d_G(v)-2)(8-d_G(v))\le9$.
Hence $|A(uv)|\ge28-(6+8+9)=5$,
and similarly $|A(uw)|\ge5$.

The list bounds $2,2,5,5$ allow an extension
by \cref{lem:hall}, a contradiction.
\end{proof}

We first prove two lemmas for extending a partial coloring.
In the first, pairs of nonconflicting edges may receive
a common color.
The second will be used for the four external edges
in the last $K_{2,3}$ case.

\begin{lemma}\label{lem:paired-clique-app}
Let $C,P_1,\ldots,P_m$ be pairwise disjoint sets
of uncolored edges
in a partial strong $28$-edge-coloring
of $G$, where $|C|=c$ and
$P_i=\{e_i,f_i\}$ for $1\le i\le m$.
Suppose that $e_i$ and $f_i$ are nonconflicting
for every $1\le i\le m$.
If $|A(e)|\ge c+m$ for $e\in C$ and
\[
 |A(e_i)|,|A(f_i)|
 \ge \left\lceil\frac c2\right\rceil+m
 \qquad (1\le i\le m),
\]
then the partial coloring extends to
$C\cup\bigcup_{i=1}^m P_i$.
\end{lemma}

\begin{proof}
Choose a maximum subfamily of these pairs that can be
colored with pairwise distinct common colors,
and let $k$ denote its size.
Color these pairs and delete their $k$ colors from
every remaining list, whether or not a particular
deletion is forced by a conflict.
Put $r=m-k$.
Each remaining list for an edge of $C$ has size
at least $c+r$, and each list for an edge in an
uncolored pair has size at least $\lceil c/2\rceil+r$.
Maximality implies that the two current lists in
each uncolored pair are disjoint.

There are $c+2r$ uncolored edges remaining.
A family containing both lists from one uncolored pair
has union of size at least
\[
 2\left(\left\lceil\frac c2\right\rceil+r\right)
 \ge c+2r.
\]
A family containing no complete pair has at most
$c+r$ members.
If it contains a list for an edge of $C$, that list
has size at least $c+r$.
Otherwise, the family has at most
$r$ members, each with a list of size at least $r$.
Thus Hall's theorem colors all remaining edges
with distinct colors, none among the $k$ deleted colors.
Together with the $k$ pairs already colored,
this gives the required extension.
\end{proof}

We will repeatedly use the following
observation.
If an uncolored edge conflicts with at most $b$ colored
edges, then it has at least $28-b$ available colors.
Uncoloring $s$ of these edges leaves at most $b-s$
colored conflicting edges and hence at least
$28-b+s$ available colors.
Since this estimate counts edges, the removed colors
need not be distinct.

\begin{lemma}\label{lem:four-tail-app}

Let $C,H,P,Q$ be pairwise disjoint sets of uncolored
edges in a partial strong $28$-edge-coloring of $G$,
where $C$ consists of six pairwise conflicting edges,
$H=\{h_1,h_2\}$, $P=\{p_1,p_2\}$, and $Q=\{q_1,q_2\}$.
Suppose that $h_1,h_2$ are nonconflicting,
each $h_i$ conflicts with every edge in $C\cup P\cup Q$,
and the two edges within each of $P$ and $Q$ conflict.
Let $J$ be the bipartite graph with parts $P,Q$,
where $p_iq_j\in E(J)$ exactly when $p_i,q_j$ conflict.
If $J\ne K_{2,2}$ and
\[
  |A(e)|\ge
  \begin{cases}
    10,       & e\in C,\\
    7,        & e\in H,\\
    6+d_J(e), & e\in P\cup Q,
  \end{cases}
\]
then the partial coloring extends to $C\cup H\cup P\cup Q$.
\end{lemma}

\begin{proof}

Suppose first that $A(h_1)\cap A(h_2)\ne\varnothing$,
and color $h_1,h_2$ alike.
The remaining lists have size at least nine on $C$
and at least five on $P\cup Q$.
Every proper subfamily satisfies Hall's condition.
If the full family also satisfies it, we are done;
otherwise its union is a $9$-element set $T$.
Choose a nonedge $pq$ of $J$.
The lists $A(p),A(q)\subseteq T$ each have size
at least five, so they intersect.
Color $p,q$ alike.
The eight remaining lists have size at least eight
on $C$ and at least four on the other two edges,
so Hall's theorem completes the coloring.

Now suppose that $A(h_1)\cap A(h_2)=\varnothing$.
A family containing both lists has at least fourteen
colors; one containing no list from $C$ and at most
one from $H$ has at most five members.
Since every list on $C$ has size at least ten,
only an $11$-member family consisting of
$C\cup P\cup Q$ and one edge of $H$ can fail
Hall's condition.
If no such failure occurs, we are done.
Otherwise, the union of a failing family is a
$10$-element set $T$.
Choose a nonedge $pq$ of $J$.
The lists $A(p),A(q)\subseteq T$ each have size
at least six, so they intersect.
Color $p,q$ alike.

The remaining lists have size at least nine on $C$,
at least six on $H$, and at least five on the other
two edges; the two lists on $H$ remain disjoint.
A family containing both has at least twelve colors.
Every other family either contains a list from $C$
and has at most nine members, or has at most three
members.
Thus Hall's theorem completes the coloring.
\end{proof}

\begin{lemma}\label{lem:k24-tight-app}
The graph $G$ contains no copy of $K_{2,4}$ whose two
vertices in the bipartition class of size two both
have degree four.
\end{lemma}

\begin{proof}

Suppose such a copy exists, with bipartition classes
$U=\{u,v\}$ and $W=\{w_1,w_2,w_3,w_4\}$.
By minimality, color $G-U$, leaving the eight
pairwise conflicting edges
$X=\{uw_i,vw_i:1\le i\le4\}$ uncolored.
Put
\[
  Z_i=\{w_ix:x\in N_G(w_i)\setminus U\},
  \qquad Z=\bigcup_{i=1}^4 Z_i.
\]
The colored edges conflicting with either $uw_i$
or $vw_i$ are precisely those in
\begin{equation}\label{eq:k24-blockers}
  F_i=Z\cup
  \bigcup_{w_ix\in Z_i}
  \bigl(E_{G-U}(x)\setminus\{w_ix\}\bigr).
\end{equation}

Let $k$ be the number of $5$-vertices in $W$.
Since $G$ is triangle-free, $W$ is independent,
so $|Z|=8+k$.
Consequently,
\[
  |F_i|\le
  \begin{cases}
    11+2\cdot4=19, & d_G(w_i)=4,\\
    12+3\cdot3=21, & d_G(w_i)=5.
  \end{cases}
\]
The first bound uses $k\le3$; the second uses the fact
that every neighbor of a $5$-vertex has degree four.
Thus every edge of $X$ has at least seven available colors.

Hall's condition can fail only for the full family,
in which case all eight lists have size seven.
Each $F_i$ then uses $21$ colors, forcing
$d_G(w_i)=5$ and $|F_i|=21$ for every $i$.
Equality also implies that the twelve edges of $Z$
have distinct endpoints outside $W$; otherwise a
term in the second union in \eqref{eq:k24-blockers}
would intersect $Z$.
In particular, each such endpoint has exactly one
neighbor in $W$.

Choose $p=w_ix\in Z_i$.
Since $d_G(x)=4$, among the nine endpoints outside $W$
of edges in $Z\setminus Z_i$ there is a vertex $y$
nonadjacent to $x$.
Let $q=w_jy$ be its edge in $Z$.
Since $W$ is independent, $x,y$ each have a unique
neighbor in $W$, and $xy\notin E(G)$,
the edges $p,q$ are nonconflicting.

Uncolor $p,q$.
Both belong to every $F_i$, so every edge of $X$
now has at least $28-(21-2)=9$ available colors.
Each of $p,q$ sees all eight uncolored edges of $X$.
By \cref{lem:strong-neighborhood},
\(
  |A(p)|,|A(q)|\ge28-(31-8)=5.
\)
Apply \cref{lem:paired-clique-app} with $c=8$, $m=1$,
$C=X$, and $P_1=\{p,q\}$.
The resulting extension contradicts the choice of $G$.
\end{proof}

\begin{lemma}\label{lem:no-bipartite-app}
The graph $G$ contains no $K_{2,3}$ or $K_{2,4}$. It also contains no $K_{2,5}$ and no $K_{3,3}$ having a bipartition class consisting of three $4$-vertices.
\end{lemma}

\begin{proof}

\emph{Case 1: $K_{3,3}$ and $K_{2,5}$.}
Suppose that a copy of $K_{3,3}$ has a bipartition class consisting of three $4$-vertices. Delete these vertices and color the remaining graph by minimality. Each of the three external edges sees all nine uncolored edges of $K_{3,3}$, so its list has size
at least $28-(31-9)=6$. For each edge of $K_{3,3}$, at most eight colored conflicts are reached through external edges at its
undeleted endpoint, four are external edges at the other two undeleted vertices, and four are reached through the external edge at its deleted endpoint. Its list therefore has size at least $28-(8+4+4)=12$. Thus \cref{lem:hall} colors all twelve deleted edges.

In a copy of $K_{2,5}$, the two vertices in the
bipartition class of size two have degree five,
and their five common neighbors have degree four.
Delete the two $5$-vertices and color the remaining
graph.
Each of the ten deleted edges conflicts with at most
eighteen colored edges: the ten external edges at
the common neighbors and at most eight further
edges reached through the two external edges at
its undeleted endpoint.
All ten lists have size at least ten, so
\cref{lem:hall} extends the coloring.

\smallskip
\emph{Case 2: $K_{2,4}$.}
The degree type $(4,4)$ on the bipartition class
of size two is excluded by \cref{lem:k24-tight-app}.
For degree type $(4,5)$, every common neighbor
has degree four.
Delete the $4$-vertex and color the remaining graph.
Each deleted edge has at most $21$ colored conflicts:
at most ten through the other three common neighbors
and the $5$-vertex, three at its own common neighbor,
and eight beyond that neighbor's two external edges.
Thus all four lists have size at least seven,
and \cref{lem:hall} applies.

Now let $\{u,v\}$ be the bipartition class of size two,
with $d_G(u)=d_G(v)=5$, and let
$W=\{w_1,w_2,w_3,w_4\}$ be the common neighbors.
Every $w_i$ has degree four.
Write $ux,vy$ for the remaining incident edges.
If $x=y$, there is an excluded $K_{2,5}$.
Hence $x\ne y$.
Delete $u,v$ and color the remaining graph.

Each of the eight edges of $K_{2,4}$ has at most
$8+3+8=19$ colored conflicts: eight external edges
at $W$, three colored edges at $x$ or $y$, and
eight further edges beyond the two external edges
at its endpoint in $W$.
Its list therefore has size at least nine.
Each of $ux,vy$ sees all eight uncolored edges of
$K_{2,4}$, so its list has size at least five.

If $xy\notin E(G)$, the edges $ux,vy$ are nonconflicting.
Apply \cref{lem:paired-clique-app} with $c=8$, $m=1$,
the eight edges of $K_{2,4}$ as $C$, and
$P_1=\{ux,vy\}$.

Suppose that $xy\in E(G)$.
Each of $ux,vy$ now sees the other nine uncolored
local edges, so its list has size at least six.
Hall's condition can fail only for the full family.
Suppose it fails.
Then all ten lists lie in a $9$-element set $T$,
and every list for an edge of $K_{2,4}$ equals $T$.

Let $Z$ be the set of the eight external edges at $W$.
For a fixed $i$, put
\[
 P_i=
 \bigcup_{t\in N_G(w_i)\setminus\{u,v\}}
 \bigl(E_{G-\{u,v\}}(t)\setminus\{w_it\}\bigr),
\]
and let
$R_x=E_{G-\{u,v\}}(x)$ and
$R_y=E_{G-\{u,v\}}(y)$.
The colored conflicts of $uw_i$ and $vw_i$ are,
respectively,
\(Z\cup P_i\cup R_x\) and 
\(Z\cup P_i\cup R_y.\)
Each union has at most $8+8+3=19$ edges and uses
all nineteen colors in $\colset\setminus T$.
Equality implies that $Z\cup P_i$ consists of
sixteen edges with distinct colors, none of which
occurs on $R_x\cup R_y$.
But $R_x,R_y$ each contain three edges and intersect
exactly in $\{xy\}$.
Thus $R_x\cup R_y$ consists of five pairwise
conflicting edges restricted to the remaining
three colors, a contradiction.

\smallskip
\emph{Case 3: $K_{2,3}$ with a $4$-vertex in the
bipartition class of size two.}
Let $\{u,v\}$ be this class, with $d_G(u)=4$,
and let $W=\{w_1,w_2,w_3\}$ be the common neighbors.
Delete $u,v$ and color the remaining graph.

If $(d_G(u),d_G(v))=(4,5)$, the external edges at
$u,v$ have at most $22,24$ colored conflicts,
respectively, while the edges of $K_{2,3}$ at
$u,v$ have at most $18,20$, respectively.
The sorted list bounds are therefore $4,4,6,8,8,8,10,10,10$,
so \cref{lem:hall} applies.

Assume that $d_G(u)=d_G(v)=4$.
Write $ux,vy$ for the external edges, let $n$ be
the number of $5$-vertices in $W$, and put
$\varepsilon_x=d_G(x)-4$ and
$\varepsilon_y=d_G(y)-4$.
Counting the external edges at $W$, the colored
edges at $x$ or $y$, and the edges beyond the
external edges at $w_i$ gives 
\[
|A(ux)|\ge7-n-\varepsilon_x,\qquad
|A(uw_i)|\ge
 \begin{cases}
   11-n-\varepsilon_x, & d_G(w_i)=4,\\
   10-n-\varepsilon_x, & d_G(w_i)=5,
 \end{cases}
 \]
 
with symmetric bounds on the $v$-side.
If $n\le2$, the sorted list bounds
$4,4,7,7,7,7,8,8$
allow an extension by \cref{lem:hall}.
Hence every $w_i$ has degree five.

Put $T_i=N_G(w_i)\setminus\{u,v\}, Z_i=\{w_it:t\in T_i\} (1\le i\le3).$
Each $T_i$ is an independent $3$-set of $4$-vertices.
We claim that some edges in distinct $Z_i$
are nonconflicting.
Suppose otherwise.
For $i\ne j$, put $r=|T_i\cap T_j|$.
Every vertex of $T_i\setminus T_j$ must be adjacent
to every vertex of $T_j\setminus T_i$.
If $r=0$, these sets span an excluded $K_{3,3}$.
If $r=1$, the two vertices of $T_i\setminus T_j$
form the smaller class of a $K_{2,3}$ of the
already excluded $n=1$ type.
If $r=3$, the vertices $w_i,w_j$ lie in an
excluded $K_{2,5}$.
Thus $|T_i\cap T_j|=2$ for all $i\ne j$.
The triple intersection has size one or two.
In either case, the three vertices of
$T_1\cup T_2\cup T_3$ outside the triple intersection
are pairwise adjacent by the preceding observation.
This contradicts triangle-freeness and proves the claim.

Choose nonconflicting edges $h_1,h_2$ in distinct
$Z_i$ and uncolor them.
Each sees all eight deleted edges, so
\(
  |A(h_1)|,|A(h_2)|\ge28-(31-8)=5.
\)
If $ux,vy$ are nonconflicting, the six edges of
$K_{2,3}$ now have lists of size at least eight,
and $ux,vy$ have lists of size at least five.
Apply \cref{lem:paired-clique-app} with $c=6$, $m=2$,
and pairs $\{ux,vy\}$ and $\{h_1,h_2\}$.

Suppose that $ux,vy$ conflict.
The absence of $K_{2,4}$ excludes
$x=y$, $uy\in E(G)$, and $vx\in E(G)$.
Hence $xy\in E(G)$.
Up to symmetry, $(d_G(x),d_G(y))$ is $(4,4)$
or $(4,5)$.
Each of $ux,vy$ sees at least nine uncolored edges:
the six edges of $K_{2,3}$, the other external edge,
and $h_1,h_2$.
Together with the preceding bounds, this gives
\[
\begin{array}{c|c|c|c}
(d_G(x),d_G(y))
 &|A(ux)|&|A(vy)|
 &\text{bounds for the six }K_{2,3}\text{-edges}\\ \hline
(4,4)&7&7&9,9,9,9,9,9\\
(4,5)&7&6&9,9,9,8,8,8
\end{array}
\]

If $A(h_1)\cap A(h_2)\ne\varnothing$, color
$h_1,h_2$ alike.
The sorted list bounds for the remaining eight
edges are $6,6,8,8,8,8,8,8$ or $5,6,7,7,7,8,8,8$,
so \cref{lem:hall} applies.
Otherwise, seek distinct colors for all ten edges.
Families of at most five lists satisfy Hall's condition;
families of sizes six through eight contain a list
of size at least eight, and every nine-member family
contains a list of size at least nine.
The full union has at least ten colors because
$A(h_1),A(h_2)$ are disjoint.
Thus Hall's theorem completes this case.

\smallskip
\emph{Case 4: $K_{2,3}$ with degree type $(5,5)$.}
Finally, assume that $d_G(u)=d_G(v)=5$.
Every common neighbor and every external neighbor
of $u$ or $v$ has degree four.
Let $P=\{p_1,p_2\}$ and $Q=\{q_1,q_2\}$ be the
external edges at $u,v$, respectively.
Let $J$ be the bipartite graph on $P,Q$ whose
edges represent conflicts.
Since $K_{2,4}$ is excluded, an edge of $P$ conflicts
with an edge of $Q$ exactly when their external
endpoints are adjacent.
Moreover, $J\ne K_{2,2}$: otherwise the two external
neighbors of $u$ would be $4$-vertices having the
three common neighbors consisting of $u$ and the
two external neighbors of $v$, contrary to Case~3.

Put  $T_i=N_G(w_i)\setminus\{u,v\}$, $Z_i=\{w_it:t\in T_i\} (1\le i\le3)$.
Each $T_i$ has size two.
These sets are pairwise disjoint, since a common
vertex would give a $K_{2,3}$ whose smaller class
consists of the $4$-vertices $w_i,w_j$.

Suppose every edge in $Z_i$ conflicts with every
edge in $Z_j$ for some $i\ne j$.
Then $T_i,T_j$ span a $K_{2,2}$.
Since no two $5$-vertices are adjacent, one of
$T_i,T_j$ consists of two $4$-vertices.
Those two vertices have three common neighbors:
their common neighbor in $W$ and the two vertices
in the other set.
This again contradicts Case~3.
Hence we can choose nonconflicting edges
$h_1,h_2$ in distinct $Z_i$.

By minimality, color $G-\{u,v\}$, and then uncolor
$h_1,h_2$.
Let $C$ be the six edges of $K_{2,3}$.
Before uncoloring $h_1,h_2$, each edge of $C$
has at most $6+6+8=20$ colored conflicts:
six external edges at $W$, six colored edges at
the two external neighbors of its endpoint in
$\{u,v\}$, and eight further edges beyond the
external edges at its endpoint in $W$.
Thus $|A(e)|\ge10$ for $e\in C$.

Each edge $e\in P\cup Q$ sees the six edges of $C$,
the other edge at the same $5$-vertex,
the $d_J(e)$ conflicting edges in the other part,
and $h_1,h_2$.
Hence it has at least $9+d_J(e)$ uncolored
strong neighbors, giving
$|A(e)|\ge28-\bigl(31-(9+d_J(e))\bigr)=6+d_J(e).$
Finally, each $h_i$ sees all ten deleted edges, so
$|A(h_i)|\ge28-(31-10)=7.$
All hypotheses of \cref{lem:four-tail-app} hold.
The resulting extension contradicts the choice of $G$,
completing the proof.
\end{proof}

\subsection{Four-cycles}

\begin{lemma}\label{lem:no-c4-app}
The graph $G$ contains no $4$-cycle.
\end{lemma}

\begin{proof}

Suppose that $C=v_1v_2v_3v_4v_1$ is a $4$-cycle,
with subscripts read modulo four.
By \cref{lem:no-triangle-app}, $C$ has no chord.
Put
\(
  U_i=N_G(v_i)\setminus V(C),
 \)
 \(
 X_i=\{v_ix:x\in U_i\}.
\)
The sets $U_i$ are pairwise disjoint:
a common neighbor of adjacent cycle vertices gives
a triangle, and a common neighbor of opposite
cycle vertices gives a $K_{2,3}$, contrary to
\cref{lem:no-bipartite-app}.

For opposite cycle vertices $v_i,v_j$, edges
$v_ix\in X_i$ and $v_jy\in X_j$ conflict exactly
when $xy\in E(G)$.
Not all such pairs conflict.
Otherwise, two vertices of $U_i$ would have the
three common neighbors consisting of $v_i$ and
any two vertices of $U_j$, giving a $K_{2,3}$.

Since no two $5$-vertices are adjacent, the degree
sequence of $C$ is, up to symmetry,
$(4,4,4,4)$, $(5,4,4,4)$, or $(5,4,5,4)$.

\smallskip
\emph{Case 1: $(4,4,4,4)$.}
Delete $v_2,v_4$ and color the remaining graph
by minimality.
For an external edge at a deleted vertex,
at most eight colored conflicts are incident with
the two external neighbors of that vertex,
twelve further conflicts are reached through
its own external endpoint, and four are incident
with $v_1,v_3$.
Its list therefore has size at least
$28-(8+12+4)=4$.

For a cycle edge, say $v_1v_2$, count the two
retained edges at $v_1$, at most eight further
edges through their other endpoints, eight edges
at the external neighbors of $v_2$, and the two
retained edges at $v_3$.
Its list has size at least $28-(2+8+8+2)=8$.
Thus the sorted list bounds are $4,4,4,4,8,8,8,8$,
and \cref{lem:hall} extends the coloring.

For the remaining cases, delete $V(C)$ and color
$G-V(C)$ by minimality.
Put
$D_i=\sum_{x\in U_i}(d_G(x)-1).$
Then $D_i\le8$ if $d_G(v_i)=4$, and $D_i=9$
if $d_G(v_i)=5$.
An external edge $v_ix$ has at most
\[
  D_i+(d_G(x)-1)(8-d_G(x))=D_i+12
\]
colored conflicts, while a cycle edge $v_iv_{i+1}$
has at most $D_i+D_{i+1}$.
Hence external-edge lists have size at least eight
at $4$-vertices and seven at $5$-vertices.
Cycle-edge lists have size at least twelve for
a $(4,4)$-edge and eleven otherwise.

\smallskip
\emph{Case 2: $(5,4,4,4)$.}
There are nine external edges and four cycle edges.
Choose nonconflicting edges $e\in X_2$ and
$f\in X_4$, as justified above.
Both lists have size at least eight.
If they intersect, color $e,f$ alike and delete
that color from every remaining list.
The sorted bounds for the eleven remaining lists are
$ 6,6,6,7,7,7,7,10,10,11,11,$
so \cref{lem:hall} applies.

Otherwise, $|A(e)\cup A(f)|\ge16$.
Families of at most eight lists satisfy Hall's
condition, since every list has size at least seven
and only three can be smaller than eight.
A family of nine through eleven lists either
contains a cycle-edge list of size at least eleven
or consists of all nine external lists and
contains both $A(e),A(f)$.
A family of twelve lists contains a cycle-edge list
of size at least twelve, and the full family
contains $A(e)\cup A(f)$.
Thus Hall's theorem extends the coloring.

\smallskip
\emph{Case 3: $(5,4,5,4)$.}
There are two disjoint nonconflicting pairs between
$X_1$ and $X_3$.
Indeed, otherwise all edges of the bipartite graph
of nonconflicting pairs would be incident with
a single vertex.
Two elements in one part would then conflict
with all three elements in the other.
The corresponding external endpoints would form
a $K_{2,3}$, a contradiction.
Also choose a nonconflicting pair between $X_2$
and $X_4$.
These three pairs are disjoint.

Choose a largest subfamily of these pairs that
can receive distinct common colors, and let its
size be $k$.
Color the chosen pairs and delete their $k$ colors
from every remaining list.
The two reduced lists of each unchosen pair
are disjoint; otherwise another pair could be added.
There remain $14-2k$ edges, with lists of size
at least $7-k$ on external edges and $11-k$
on cycle edges.

Let $\mathcal L$ be a nonempty family of remaining lists.
If it contains both lists of an unchosen pair,
their union has size at least $2(7-k)=14-2k$.
If it contains no complete pair but contains
a cycle-edge list, then
$ |\mathcal L| \le(14-2k)-(3-k)=11-k,$
which is at most the size of that list.
If it contains neither a complete pair nor
a cycle-edge list, then
$ |\mathcal L| \le(10-2k)-(3-k)=7-k,$
which is at most the size of each list.
Hall's theorem colors the remaining edges with
distinct colors from their reduced lists,
a contradiction.
\end{proof}

\subsection{Five-cycles}

\begin{lemma}\label{lem:no-c5-app}
The graph $G$ contains no $5$-cycle.
\end{lemma}

\begin{proof}

Suppose that $C=v_1v_2v_3v_4v_5v_1$ is a $5$-cycle,
with subscripts read modulo five.
Since no two $5$-vertices are adjacent, we may
relabel $C$ so that
$d_G(v_2)=d_G(v_4)=d_G(v_5)=4$.
Put $c_i=v_iv_{i+1}$.
By \cref{lem:no-triangle-app,lem:no-c4-app},
$C$ is chordless and no vertex outside $C$
has two neighbors on $C$.

By minimality, fix a strong $28$-edge-coloring
$\varphi$ of $G-V(C)$.
Since this subgraph is induced, $\varphi$ is
a partial strong edge-coloring of $G$.
Write $A(e)=A_\varphi(e)$ for each uncolored edge $e$.

For an external edge $v_ix$, the colored conflicts
consist of edges incident with the external
neighbors of $v_i$ and further edges reached
through $x$.
The first group has size at most eight if
$d_G(v_i)=4$ and nine otherwise.
The second has size at most
$(d_G(x)-1)(8-d_G(x))=12$.
A cycle edge has at most seventeen colored conflicts.
Hence
\begin{equation}\label{eq:c5-repair-lists}
 |A(v_ix)|\ge
 \begin{cases}
   8, & d_G(v_i)=4,\\
   7, & d_G(v_i)=5,
 \end{cases}
 \qquad
 |A(c_i)|\ge11.
\end{equation}

Every external edge at $v_i$ is nonconflicting
with $c_{i+2}$.
Between the external edges at two nonadjacent
cycle vertices, the conflict relation is a matching:
two such edges conflict exactly when their outside
endpoints are adjacent, and an outside endpoint
adjacent to two outside neighbors of one cycle
vertex would create a $4$-cycle.

If $v_1$ or $v_3$ has degree four, supplement its
two external edges with one formal element,
having an arbitrary seven-element list in $\colset$
and declared nonconflicting with every other element.
These formal elements do not alter $G$.
Label the resulting elements $e_{ij}$, with three
at each of $v_1,v_3$ and two at each other cycle vertex.
Keep all lists $A(e)$ fixed throughout.

By the matching property, choose the labels so that
conflicts between the groups at $v_1,v_3$,
and between those at $v_2,v_4$, occur only between
elements with the same second subscript.
Define
\begin{equation}\label{eq:c5-repair-pairs}
\begin{aligned}
 \mathcal P
 &=\bigl\{\{e_{11},e_{32}\},\{e_{12},e_{33}\},
          \{e_{21},e_{42}\},\{e_{22},e_{41}\}\bigr\},\\
 \mathcal Q
 &=\bigl\{\{e_{13},c_3\},\{e_{31},c_5\},
          \{e_{51},c_2\}\bigr\}.
\end{aligned}
\end{equation}
These seven pairs are disjoint and have
nonconflicting members.
The three unpaired elements are $e_{52},c_1,c_4$.
Every external edge or formal element has a list
of size at least seven, and the actual edge
$e_{52}$ has at least eight.

Choose as many of these pairs as possible so that
each receives a common available color and
different pairs receive distinct colors.
Subject to this, maximize the number chosen
from $\mathcal P$.
Let $k$ be the total number chosen and $h$ the
number chosen from $\mathcal Q$.
Assign these colors to the chosen pairs.
Obtain $L(e)$ from $A(e)$ by deleting all $k$
used colors.
There remain $17-2k$ elements.
By maximality, the two reduced lists of every
unchosen pair are disjoint.

For an unchosen pair $\{e,f\}\in\mathcal P$,
no color used on a chosen pair of $\mathcal Q$
belongs to $A(e)\cap A(f)$.
Otherwise, replacing that chosen pair by $\{e,f\}$
would preserve $k$ and increase the number chosen
from $\mathcal P$.
Thus the $h$ colors used on $\mathcal Q$ reduce
the total size of these two lists by at most $h$,
and the other colors reduce it by at most $2(k-h)$.
Consequently,
\begin{equation}\label{eq:c5-hall-priority}
 |L(e)\cup L(f)|
 =|L(e)|+|L(f)|
 \ge14-2k+h
 \qquad
 (\{e,f\}\in\mathcal P\text{ unchosen}).
\end{equation}
For an unchosen pair in $\mathcal Q$, the union
has size at least $7+11-2k=18-2k$.

Let $Z$ be a nonempty set of remaining elements.
If $Z$ contains a complete pair of $\mathcal Q$,
its list union has size at least
$18-2k>17-2k\ge|Z|$.
If it contains a complete pair of $\mathcal P$
but none of $\mathcal Q$, it omits at least one
element from each of the $3-h$ unchosen pairs
of $\mathcal Q$.
Hence
$|Z|\le17-2k-(3-h)=14-2k+h$,
and \eqref{eq:c5-hall-priority} gives Hall's condition.

Finally, suppose that $Z$ contains no complete pair.
Then $|Z|\le(7-k)+3=10-k$.
If $Z$ contains a cycle edge, its list has size
at least $11-k$, which suffices.
Otherwise, the only unpaired element that can
belong to $Z$ is $e_{52}$.
If $e_{52}\in Z$, then
$|Z|\le8-k\le|L(e_{52})|$.
If not, then $|Z|\le7-k$, and every remaining
list has size at least $7-k$.
Thus Hall's condition holds in all cases.

By Hall's theorem, the remaining elements receive
distinct colors avoiding those used on the chosen
pairs.
Discarding the formal elements gives a strong
$28$-edge-coloring of $G$ extending $\varphi$,
a contradiction.
\end{proof}

Together, \cref{lem:no-triangle-app,lem:no-c4-app,lem:no-c5-app}
establish \cref{cl:girth}.

\section*{Data availability}
No datasets were generated or analyzed in this study.

\section*{Declaration of competing interest}
The authors declare that they have no competing interests.

\end{document}